\documentclass{amsart}
\usepackage{graphicx}
\usepackage{amssymb}
\usepackage{epstopdf}
\usepackage{nicefrac}
\usepackage{enumerate}
\usepackage{float}
\usepackage{amsaddr}
\usepackage{hyperref}

\newtheorem{thm}{THEOREM}[section]

\newtheorem{cor}[thm]{COROLLARY}
\newtheorem{defn}[thm]{DEFINITION}

\newtheorem{lemma}[thm]{LEMMA}
\newtheorem{prob}[thm]{PROBLEM}
\newtheorem{prop}[thm]{PROPOSITION}

\newtheorem{remark}[thm]{REMARK}

\newcommand{\G}{\Gamma}

\newcommand{\cF}{{\mathcal F}}
\newcommand{\cR}{{\mathcal R}}

\newcommand{\wth}{{\widetilde{h}}}

\newcommand{\mG}{{\mathbb G}}

\newcommand{\mN}{{\mathbb N}}

\newcommand{\mR}{{\mathbb R}}

\newcommand{\mZ}{{\mathbb Z}}

\newcommand{\cA}{{\mathcal A}}

\newcommand{\cC}{{\mathcal C}}

\newcommand{\cE}{{\mathcal E}}
\newcommand{\cG}{{\mathcal G}}

\newcommand{\cM}{{\mathcal M}}

\newcommand{\cO}{{\mathcal O}}
\newcommand{\cP}{{\mathcal P}}
\newcommand{\cQ}{{\mathcal Q}}
\newcommand{\cS}{{\mathcal S}}
\newcommand{\cT}{{\mathcal T}}
\newcommand{\cU}{{\mathcal U}}
\newcommand{\cV}{{\mathcal V}}

\newcommand{\cX}{{\mathcal X}}
\newcommand{\cZ}{{\mathcal Z}}

\newcommand{\fG}{\mathfrak{G}}

\newcommand{\fM}{{\mathfrak{M}}}

\newcommand{\fP}{{\mathfrak{P}}}

\newcommand{\fX}{{\mathfrak{X}}}

\input xy
\xyoption {all}

\begin{document}

\title{Hierarchy of graph matchbox manifolds}

\thanks{2000 {\it Mathematics Subject Classification}. Primary 57R30. Secondary 37B05, 54H20}

\author{Olga Lukina}

\address{Department of Mathematics, University of Leicester, University Road, \\ Leicester LE1 7RH, United Kingdom}
\email{ollukina940@gmail.com}

\thanks{Version date: August 9, 2012}

\date{}

\keywords{matchbox manifold, laminations, ends of leaves, growth of leaves, theory of levels, Gromov-Hausdorff metric}

\begin{abstract}
We study a class of graph foliated spaces, or graph matchbox manifolds, initially constructed by Kenyon and Ghys. For graph foliated spaces we introduce a quantifier of dynamical complexity which we call its level. We develop the \emph{fusion} construction, which allows us to associate to every two graph foliated spaces a third one which contains the former two in its closure. Although the underlying idea of the fusion is simple, it gives us a powerful tool to study graph foliated spaces. Using fusion, we prove that there is a hierarchy of graph foliated spaces at infinite levels. We also construct examples of graph foliated spaces with various dynamical and geometric properties. 
\end{abstract}

\maketitle

\section{Introduction} \label{sec-intro}

A \emph{matchbox manifold} is a compact connected metrizable space $M$ such that each point $x \in M$ has a neighborhood homeomorphic to a product space $U_x \times N_x$, where $U_x \subset \mR^n$ is open and $N_x$ is a compact totally disconnected space. The term `matchbox manifold' originates from the works of Aarts and Martens \cite{AM1988}, Aarts and Oversteegen \cite{AO1995} for the case when $n=1$, when local charts can be thought of as `boxes of matches'. The most well-studied classes of examples of matchbox manifolds are weak solenoids \cite{McCord1965,FO2002}, generalized solenoids \cite{Williams1974}, and tiling spaces of aperiodic tilings with finite local complexity (see, for instance, \cite{FrankSadun2009}, or \cite{BG2003} for a more general type of tilings). In this paper we consider a third  class of examples, which we call \emph{graph matchbox manifolds}. This construction was introduced by  Kenyon and Ghys \cite{Ghys1999}, and later generalized by Blanc \cite{Blanc2001}, Lozano Rojo \cite{LozanoSp}, Alcalde Cuesta, Lozano Rojo and Macho Stadler \cite{ALM2008}. 

\begin{remark}[On the use of terminology]
{\rm 
The notion of a matchbox manifold is essentially the same as that of a \emph{lamination}. The term `lamination' appears in the literature in two slightly different contexts: in low-dimensional topology, a lamination is a decomposition into leaves of a closed subset of a manifold; in holomorphic dynamics, Sullivan \cite{Sullivan1988} introduced Riemann surface laminations as compact topological spaces locally homeomorphic to a complex disk times a Cantor set. An embedding into a manifold is not required in the latter context, and a matchbox manifold is a lamination in this terminology. The concept of a \emph{foliated space} as a generalization of a foliated manifold was introduced in the book by Moore and Schochet \cite{MS2006}, where a foliated space is defined as a separable metrizable space locally homeomorphic to a product of a disk in $\mathbb{R}^n$ and a separable metrizable space. In this terminology, a matchbox manifold is a foliated space with specific properties, i.e. it is a compact foliated space with totally disconnected transversals. In the present paper we follow terminology of Candel and Conlon \cite{CandelConlon2000}, reserving the word `lamination' for a foliated space embedded in a manifold. We then use the term `matchbox manifold' to distinguish a class of foliated spaces which are compact and have totally disconnected transversals.
}
\end{remark}

Let $G$ be a finitely generated group with a non-symmetric set of generators $G_0$, that  is, if $h \in G_0$ then $h^{-1} \notin G_0$. Let $\cG$ be the Cayley graph of $G$, and $X$ be the set of all infinite connected subtrees of $\cG$ containing the identity $e$. Each subtree $T$ is equipped with a standard complete length metric $d$, and the pair $(T,e)$ is a pointed metric space. The set $X$, endowed with the Gromov-Hausdorff metric $d_{GH}$ \cite{Burago}, is a compact totally disconnected space \cite{Blanc2001,Ghys1999,LozanoSp}. One can define a partial action of the free group $F_n$ on $X$, where $n$ is the cardinality of the set of generators $G_0$. This action gives rise to a pseudogroup $\fG$ on $X$, and an important feature of the construction is that the pseudogroup dynamical system $(X,\fG)$ can be realised as the holonomy pseudogroup of a smooth foliated space $\fM_G$ with $2$-dimensional leaves \cite{Blanc2001,Ghys1999,LozanoSp}. By this construction, for $(T,e) \in X$ the corresponding leaf $L_T \subset \fM_G$ can be thought of as the two-dimensional boundary of the thickening of a quotient graph of $T$, where the quotient map is determined by the geometry of $T$. 

\begin{defn}\label{defn-GMM}
A \emph{graph matchbox manifold} is the closure $\cM = \overline{L}$ of a leaf $L$ in $\fM_G$, that is, $\cM$ is a closed saturated transitive subset of $\fM_G$.
\end{defn}

\subsection{Hierarchies of graph matchbox manifolds} In previous works the construction of Kenyon and Ghys was mostly used to produce examples of matchbox manifolds with specific geometric and ergodic properties. Ghys \cite{Ghys1999}, see also \cite{ALM2008}, showed that if $G=\mZ^2$ then $\fM_{\mZ^2}$ contains a leaf $L$ such that the matchbox manifold $\cM=\overline{L}$ is minimal and has leaves with different conformal structures. Lozano Rojo \cite{Lozano2009} studied minimal examples in the case $G = \mZ^2$ from the point of view of ergodic theory. In the case where $G = F_3$, a free group on three generators, Blanc \cite{Blanc2001} found an example of a graph matchbox manifold containing leaves with any possible number of ends. 

In this paper we study a partial order on the corresponding foliated space $\fM_G$, given by inclusions. The following basic observation allows to restrict our attention to the case $G =F_n$ and the space of graph matchbox manifolds $\fM_n$. 

\begin{thm}\label{prop-univprop}
Given a group $G$ with a set of generators $G_0$ of cardinality at most $n$, there exists a foliated embedding 
  $$\Phi:\fM_G \to \fM_n,$$
where $\fM_G$ and $\fM_n$ are foliated spaces obtained by the construction of Kenyon and Ghys for $G$ and a free group $F_n$ on $n$ generators respectively.
\end{thm}

Let $\cM_1,\cM_2 \subset \fM_n$ be graph matchbox manifolds, then the rule 
 \begin{align*}\cM_1 \preceq \cM_2 ~ \textrm{ if and only if}~ \cM_1 \subseteq \cM_2\end{align*}
defines a partial order on the set $\cS_n$ of graph matchbox manifolds in $\fM_n$. Compact leaves and minimal subsets of $\fM_n$ are minimal elements in $\cS_n$ with respect to this order. The following theorem describes the structure of $\fM_n$.

Recall \cite{Blanc2001} that a leaf $L \subset \cM$ is \emph{recurrent} if and only if $L$ is transitive and accumulates on itself. A leaf $L$ is \emph{proper} if it does not accumulate on itself.

\begin{thm}\label{thm-globalleaf}
The partially ordered set $(\cS_n,\preceq)$ of graph matchbox manifolds in the foliated space $\fM_n$, $n>1$, has the following properties.
\begin{enumerate}
\item the set $C = \{L \subset \fM_n ~|~ L ~ \textrm{is compact}\}$ is a dense meager subset of $\fM_n$. Moreover, $C \cap \cX$ is countable, where $\cX$ is a canonical embedding of $X$ into $\fM_n$.
\item $(\cS_n, \preceq)$ is a directed partially ordered set, i.e. given $\cM_1, \cM_2 \in \cS_n$ there exists $\cM_3 \in \cS_n$ such that $\cM_1 \cup \cM_2 \subseteq \cM_3$.
\item $(\cS_n, \preceq)$ contains a unique maximal element $\cM_{max} = \fM_n$ which has a recurrent leaf. Therefore, $\fM_n$ contains a residual subset of recurrent leaves.
\end{enumerate}
\end{thm}

In order to prove Theorem \ref{thm-globalleaf}.(2), we introduce the `fusion' construction which associates to any two transitive subsets $\cM_1$ and $\cM_2$ of $\fM_n$ a transitive subset $\cM_3$ such that $\cM_3 \supseteq \cM_1 \cup \cM_2$. More precisely, given pointed graphs $(T_1,e)$ and $(T_2,e)$ such that $\cM_1 = \overline{L_{T_1}}$ and $\cM_2 = \overline{L_{T_2}}$ we give a recipe to construct a graph $(T_3,e)$ such that $\cM_3 = \overline{L_{T_3}}$ satisfies the required property. The underlying idea of the construction is very simple, but it gives us a powerful tool which allows us to obtain a lot of information about hierarchy and properties of graph matchbox manifolds. Theorems \ref{thm-globalleaf}.(2) and \ref{thm-globalleaf}.(3) are the first applications of fusion. 

Theorem \ref{thm-globalleaf}.(2) is a direct consequence of the fusion. The next important observation is that in the space $\fM_n$ fusion enables us to construct infinite increasing chains of graph matchbox manifolds. Using \cite[Theorem 3.5]{Blanc2001} with slightly eased assumptions, we conclude that the closure of such a chain contains a dense leaf, and if every element in the chain is distinct, this dense leaf accumulates on itself. The existence of a maximal closed recurrent subset in the space of graph matchbox manifolds follows by standard topological arguments, and, using fusion again, one argues that if a maximal recurrent subset exists, then it is unique, which is the statement of Theorem \ref{thm-globalleaf}.(3). 

Our main theorem shows that the space $\fM_n$ of graph matchbox manifolds contains a complicated hierarchy of finite and infinite chains of distinct graph matchbox manifolds, which motivates further the study of \emph{theory of levels} for graph matchbox manifolds.

\begin{thm}\label{thm-infiniteleaves} For a space of graph matchbox manifolds $\fM_n$, $n >1$, the following holds.
\begin{enumerate}
\item There exists an infinite increasing chain
  \begin{align*} \cM_0 \subset \cM_1 \subset \cM_2 \subset \cdots\end{align*}
of distinct graph matchbox manifolds such that $\cM = \overline{\bigcup_i \cM_i}$ is a proper subset of $\fM_n$. 
\item Let $\cM$ be a graph matchbox manifold, and suppose $\cM$ is a proper subset of $\fM_n$. Then there exists a graph matchbox manifold $\widetilde{\cM}$ such that 
  $$\cM \subset \widetilde{\cM} \subset \fM_n$$ 
are proper inclusions.
\end{enumerate}
\end{thm}

The fusion technique described earlier lies at the heart of the proof of Theorem \ref{thm-infiniteleaves}. 

In the rest of the paper we make precise the hierarchy of graph matchbox manifolds, suggested by Theorem \ref{thm-infiniteleaves}, associating to them a quantifier of dynamical complexity, called \emph{level}, inspired by the ideas of Cantwell and Conlon \cite{CC1, CC3}, Hector \cite{Hector,Hector2,Hector3}, Nishimori \cite{Nish1977,Nish1979} and Tsuchiya \cite{Ts1980}. Our main tool is the fusion technique of Theorem \ref{thm-globalleaf}.(2). Thus fusion provides the means to study the hierarchy of graph matchbox manifolds in $\fM_n$ given by inclusions. We also investigate how large is the set of matchbox manifolds realisable as graph matchbox manifolds.

\subsection{Theory of levels} The study of the partial ordering given by inclusions of transitive subsets is a natural problem  in topological dynamics. The consideration of the extension of the Poincar\'{e} Recurrence Theorem for flows to the closures of leaves of foliations began in the 1950's.
In foliation theory, the strongest results have been obtained for codimension $1$ transversally $C^2$-differentiable foliations by Cantwell and Conlon \cite{CC1, CC3}, Hector \cite{Hector, Hector2}, Nishimori \cite{Nish1977,Nish1979} and Tsuchiya \cite{Ts1980}. Later, the extent to which these ideas carry on to the codimension $1$ $C^0$ case was studied by Salhi \cite{Salhi1982,Salhi1985-1,Salhi1985-2}, and, for foliations of higher codimensions satisfying certain additional conditions, by Marzougui and Salhi \cite{Salhi2003}.
Cantwell and Conlon \cite{CC1} introduced the notion of a \emph{level} of a leaf or a transitive subset, which can be seen as a quantifier of dynamical complexity. A similar, but not the same notion of `depth' of a leaf was considered by Nishimori \cite{Nish1977,Nish1979}. The theory of levels for $C^2$-foliations relies heavily on the Kopell lemma \cite{CC1,CandelConlon2000}, which does not apply in our setting. However, it is possible to introduce the \emph{level of a graph matchbox manifold}, which is similar to the notion of a level in \cite{CC1}.

\begin{defn}\label{defn-level}
Let $\cM \subset \fM_n$ be a graph matchbox manifold. 
\begin{enumerate}
\item $\cM$ is said to be \emph{at level $0$} if either $\cM$ is a compact leaf, or $\cM$ is a minimal foliated space. In that case all leaves of $\cM$ are also at level $0$. 
\item $\cM$ is \emph{at level $k$} if the closure of the union of leaves which are not dense in $\cM$, is a proper closed subset of $\cM$, every such leaf is at level at most $k-1$, and there is at least one leaf at level $k-1$. A leaf $L$ is at level $k$ if it is dense in $\cM$.
\item $\cM$ is at \emph{infinite level} if it is not at finite level. Then a leaf $L \subset \cM$ is at infinite level if it is dense in $\cM$.
\end{enumerate}
\end{defn}

By Theorem \ref{thm-globalleaf}.(3) the space $\fM_n$ contains a dense leaf $L$ and so $\fM_n$ is a graph matchbox manifold. Since by Theorem \ref{thm-globalleaf}.(1) the set of compact leaves is dense in $\fM_n$, $\fM_n$ cannot be at  finite level, and so is at infinite level, and every dense leaf in $\fM_n$ is at infinite level. By Theorem \ref{thm-globalleaf}.(3) $\fM_n$ contains a residual subset of dense leaves, therefore, we have the following corollary.

\begin{cor}\label{thm-examples1inf}
In a foliated space $\fM_n$, $n>1$, leaves at infinite level form a residual subset.
\end{cor}

Leaves at infinite level in Corollary \ref{thm-examples1inf} are leaves which are dense in $\fM_n$, $n>1$, and so contain any other leaf in their closure. Theorem \ref{thm-infiniteleaves} shows that the meager subset of leaves which are not dense in $\fM_n$ also contains leaves at infinite level, and we are now interested in those.

Definition \ref{defn-level} allows for two types of matchbox manifolds at infinite level which are properly included in $\fM_n$: we say that a matchbox manifold $\cM$ is at infinite level of \emph{Type 1} if the union of leaves which are not dense in $\cM$ is dense in $\cM$, and $\cM$ is a proper subset of $\fM_n$. We say that a matchbox manifold $\cM$ is at infinite level of \emph{Type 2} if the union of leaves which are not dense in $\cM$ is a proper subset of $\cM$ and contains a leaf at infinite level. Our main theorem \ref{thm-infiniteleaves} states that graph matchbox manifolds at infinite level of both types exist. 

Then natural question to ask is the following. Given $k\geq 0$, does there exist a graph matchbox manifold $\cM^{(k)}$ at level precisely $k$? We now proceed to investigate these questions.

For $k=0,1$ examples of $\cM^{(0)}$ and $\cM^{(1)}$ are fairly easy to construct. Theorems \ref{thm-examples1totprop} and \ref{thm-examples1recurrent} give examples of graph matchbox manifolds at level $2$.

Recall that a leaf $L$ is proper if it does not accumulate on itself. A leaf $L$ is totally proper if $\overline{L}$ contains only proper leaves. Totally proper leaves are the simplest dynamically; the following theorem shows that in the space of graph matchbox manifolds even totally proper leaves at level $2$ can be fairly complicated, for example, they can contain an infinite number of distinct compact leaves in their closure.

\begin{thm}\label{thm-examples1totprop}
In a foliated space $\fM_n$, $n>1$, there is a totally proper graph matchbox manifold $\cM$ at level $2$ which contains a countably infinite number of compact leaves.
\end{thm}

The next theorem shows that leaves at finite level need not be totally proper, giving an example of a graph matchbox manifold at level $2$ with recurrent leaves.

\begin{thm}\label{thm-examples1recurrent}
In a foliated space $\fM_n$, there exists a graph matchbox manifold $\cM$ which contains an uncountable number of leaves, all but a finite number of which are at level $2$. The complete transversal $\fX$ of $\cM$ contains a clopen subset $\cO$ such that the restricted pseudogroup $\fG|_{\cO}$ is equicontinuous. Every leaf in $\cM$ which is at level $2$ is recurrent.
\end{thm}

Already for $k=2$, the combinatorial arguments of Theorem \ref{thm-examples1totprop} become quite tedious; it is therefore reasonable to look for topological and geometric obstructions which would allow us to determine a level of a graph matchbox manifold by other methods. 

For codimension $1$ $C^2$ foliations the growth function of a totally proper leaf at level $k$ is a polynomial of degree $k$ \cite{CC1,Ts1980} and, if the growth function of a leaf is dominated by a polynomial, then the leaf is totally proper \cite{CC1,CC3}. The situation is different in the case of graph matchbox manifolds, as the following theorem shows.

\begin{thm}\label{thm-examples1growth} 
In a foliated space $\fM_n$, there exists a totally proper graph matchbox manifold $\cM_1$ at level $1$ with a transitive leaf $L$ with linear growth, and there is also a totally proper graph matchbox manifold $\cM_2$ at level $1$ with a transitive leaf $L'$ with exponential growth. In addition, $\fM_n$ contains a recurrent graph matchbox manifold $\cM_3$ with leaves of polynomial growth.
\end{thm}

Theorem \ref{thm-examples1growth} shows that exponential growth of a dense leaf in a graph matchbox manifold cannot serve as an obstruction to the graph matchbox manifold being at finite level.

\subsection{Transverse dynamics of graph matchbox manifolds} The remaining part of the article is devoted to understanding how large is the class of examples which can be obtained by the construction of Kenyon and Ghys. We first notice that a matchbox manifold with more than one leaf must have expansive transverse dynamics, as the following proposition shows.

Let $\cO$ be a clopen subset of a transversal space $\fX$ where the Gromov-Hausdorff topology on $\fX$ is realised by the ball metric $d_X$ defined in Section \ref{subsec-kenyongraphs}. Recall \cite{Hurder2010} that the pseudogroup $\fG|_{\cO}$ is $\epsilon$\emph{-expansive} if there exists $\epsilon>0$ so that for all $w \ne w' \in \cO$ with $d_X(w,w') < \epsilon$ there exists a holonomy homeomorphism $h \in \fG|_{\cO}$ with $w,w' \in {\rm dom} (h)$ such that $d_{X}(h(w),h(w')) \geq \epsilon$. A pseudogroup $\fG|_{\cO}$ is \emph{equicontinuous} if for every $\epsilon>0$ there exists $\delta>0$ such that for all $h \in \fG|_{\cO}$ and  all $w,w' \in {\rm dom}(h)$ such that $d_{X}(w,w') < \delta$ we have $d_X(h(w),h(w'))< \epsilon$.

\begin{thm}\label{thm-dynamicsvsrecurrency}
Let $\cM$ be a graph matchbox manifold in a foliated space $\fM_n$, and $(\fX,d_X)$ be a complete transversal for $\cM$. Then the following holds.
\begin{enumerate}
\item $\cM$ contains a non-compact leaf if and only if for all $\epsilon < e^{-2}$ the restriction of the pseudogroup $\fG$ to $\fX$ is $\epsilon$-expansive. 
\item If there exists a clopen $\cO \subset \fX$ such that $\fG|_\cO$ is equicontinuous, then $\cM$ contains a recurrent leaf. The converse is false, since the restriction of the pseudogroup of any minimal matchbox manifold to any clopen subset of the transversal is expansive.
\end{enumerate}
\end{thm}
 
The result of Theorem \ref{thm-dynamicsvsrecurrency} is not surprising, since the ball metric on $X_n$ is similar to the metric on tiling spaces, which is well-studied. In particular, it is known that for a non-periodic tiling, the transverse dynamics of its tiling space is expansive (see, for example, Benedetti and Gambaudo \cite{BG2003}). Leaves in tiling spaces are usually copies of the Euclidean space, and it is customary to consider tiling spaces where each leaf is dense. To the best of our knowledge, the most general type of tiling spaces which has been described in the literature is that in Benedetti and Gambaudo \cite{BG2003}, where leaves are homogeneous spaces and are homeomorphic. Leaves in graph matchbox manifolds exhibit a much wider variety of topological properties, in particular, one can realise a surface of any genus as a leaf in $\fM_n$. Theorem \ref{thm-dynamicsvsrecurrency}.(1) uses similar arguments to those in tilings to investigate transverse dynamics for graph matchbox manifolds.

Theorem \ref{thm-dynamicsvsrecurrency}.(1) allows us to restrict the range of spaces which can be modeled by the method of Kenyon and Ghys, as the following corollary shows.

\begin{cor}\label{no-solenoid}
The foliated space $\fM_n$ does not contain a weak solenoid.
\end{cor}

Indeed, a weak solenoid is a minimal space with equicontinuous transverse dynamics. If a graph matchbox manifold $\cM$ is at level $k=0$, then it is either minimal and has expansive dynamics, or it is a compact leaf and has trivially equicontinuous dynamics. As the example in Theorem \ref{thm-examples1recurrent} shows, if $\cM$ is at level $k>0$, various types of dynamics can mix. By Theorem \ref{thm-dynamicsvsrecurrency}.(1), the clopen set $\cO$ in Theorem \ref{thm-examples1recurrent} is necessarily not a complete transversal.

\subsection{Open problems} A natural question and an interesting open problem is whether a given foliated space can be embedded as a closed saturated subset of a smooth foliated manifold. With this question in mind, it is interesting to compare the result of Theorem \ref{thm-infiniteleaves} with the hierarchy of leaves at infinite levels for codimension $1$ foliations. Namely, \cite{CC1,Hector} state that if the transverse differentiability of a codimension $1$ foliation is $C^0$, then there exists a hierarchy of infinite levels, as in Theorem \ref{thm-infiniteleaves}, but if the transverse differentiability is $C^2$, there is only one infinite level. Therefore, we conclude that graph matchbox manifolds containing leaves at more than one infinite level cannot embed as foliated subsets of codimension $1$ smooth foliations. Then the following question is natural.

\begin{prob}\label{questEmb}
Let $\cM \subset \fM_n$ be a graph matchbox manifold which contains leaves at more than one infinite levels. Does there exist a foliated embedding $\cM \to M$ into a smooth foliated manifold $M$ with a foliation of codimension $q \geq 2$?
\end{prob}

It was shown in \cite{Lozano2009} that a minimal graph matchbox manifold can always be embedded as a codimension $2$ subset of a manifold topologically. An embedding of $\cM$ in a smooth way and, moreover, as a subset of a smooth foliation, is an infinitely more subtle and technically demanding procedure. In general, not much is known about embeddings of foliated spaces transversely modeled on Cantor sets as saturated subsets of smooth foliations and, to the best of our knowledge, there are only a few works devoted to that question. For example, Williams \cite{Williams1974} showed that the expansive dynamical system associated to a transverse section of a generalized solenoid is conjugate to the dynamics of an expanding attractor, but it is an open question whether such a solenoid embeds as a subset of the corresponding foliation of a manifold. Clark and Hurder \cite{ClarkHurder2010a} give sufficient conditions under which certain homogeneous solenoids embed as subsets of smooth foliations. Problem \ref{questEmb} asks a similar question for graph matchbox manifolds. Corollary \ref{no-solenoid} shows that the space of graph matchbox manifolds does not contain weak solenoids, therefore, the result in \cite{ClarkHurder2010a} does not provide even a partial answer to Problem \ref{questEmb}. 

Following \cite{CC1} we denote by $S(L)$ the union of non-dense leaves in a graph matchbox manifold $\cM = \overline{L}$, and call $S(L)$ the substructure of a leaf $L$. In codimension $1$ $C^2$ foliations leaves at infinite levels have a substructure of leaves at finite levels, and each of the finite levels is represented in $S(L)$ by at least one leaf \cite{CC1}. In this paper we have not given an algorithm which, for a given graph matchbox manifold, determines its level. Therefore, the following questions remain.

\begin{prob}\label{various-levels}
\begin{enumerate}
\item Does there exist a graph matchbox manifold $\cM = \overline{L} \subset \fM$ at infinite level such that $S(L)$ has a substructure of leaves at all finite levels?
\item Given $k>2$, give an example of a graph matchbox manifold $\cM_k$ at level precisely $k$.
\end{enumerate}
\end{prob}

Corollary \ref{no-solenoid} shows that a large subclass of matchbox manifolds with equicontinuous dynamics, namely weak solenoids, cannot arise as graph matchbox manifolds. Theorem \ref{thm-examples1recurrent} shows that one can construct a non-minimal graph matchbox manifold which admits a transverse non-global section with equicontinuous restricted dynamics. Then there is the following question.

\begin{prob}\label{solenoid-inside}
What is precisely the relation between weak solenoids and graph matchbox manifolds?
\end{prob}

\subsection{Organization of the paper}
The rest of the paper is organized as follows. In Section \ref{subsec-kenyongraphs} we give details of the generalized construction of Ghys and Kenyon, and describe some relations between topology of trees in $X$ and the graphs of their orbits under the pseudogroup $\fG$. We also give examples of graph matchbox manifolds, and of constructing graph matchbox manifolds with desired properties, and prove Theorem \ref{prop-univprop}. We prove Theorem \ref{thm-globalleaf} in Section \ref{asymptotic}, and Theorems \ref{thm-infiniteleaves}, \ref{thm-examples1totprop} and \ref{thm-examples1recurrent} in Section \ref{exmpl}. Theorem \ref{thm-examples1growth} and Theorem \ref{thm-dynamicsvsrecurrency} are proved in Section \ref{further}. 

{\bf Acknowledgements:} The author thanks Alex Clark and Steve Hurder for their encouragements to pursue this research topic and for offering many useful comments. The author thanks the referee for recommending improvements in the presentation of this paper.

\section{Construction of a foliated space after Kenyon and Ghys}\label{subsec-kenyongraphs}

In this section we give an outline of the construction of Ghys \cite{ALM2008,Ghys1999,LozanoSp}, and define  graph matchbox manifolds.
 
 We also obtain Theorem \ref{prop-univprop}, that is, that the foliated space $\fM_n$ obtained by the construction of Kenyon and Ghys in the case $G = F_n$ has a universal property.

\subsection{Preliminaries and notation} Given a directed graph (or digraph) $T$, denote by $V(T)$ the set of vertices of $T$, and by $E(T)$ the set of edges of $T$. For an edge $w \in E(T)$ denote by $s(w)$ and $t(w)$ respectively its starting and its ending vertex. A subgraph of $T$ is a graph $T'$ with the set of vertices $V(T') \subseteq V(T)$ and the set of edges $E(T') \subset E(T)$, where an edge $w \in E(T')$ if and only if $s(w) \in V(T')$ and $t(w) \in V(T')$. A labeling of a graph $T$, or of the set $E(T)$ of edges of $T$, by a set $\cA$, is given by a function $a: E(T) \to \cA$. If it is necessary to keep track of labeling of edges in $E(T)$, we use the notation $w_{a(w)}$ for an edge $w \in E(T)$.

Let $G$ be a finitely generated group acting on itself on the right, and let $G_0$ denote a set of generators of $G$. We assume that $G_0$ is not a symmetric set, that is, if $h \in G_0$ then $h^{-1} \notin G_0$. Let $\cG$ be the Cayley graph of $G$ with the set of generators $G_0$; more precisely, set $V(\cG) = G$, and to each pair $g_1,g_2 \in V(\cG)$ such that $g_1 h = g_2$ for some $h \in G_0$, associate an edge $w_h \in E(\cG)$ with $s(w_h) = g_1$ and $t(w_h) = g_2$. Thus $\cG$ is a directed graph labeled by the set $G_0$. 

We define a length structure $\ell: P(\cG) \to \mR: \delta \mapsto \ell(\delta)$, where $P(\cG)$ the set of all paths in $\cG$, in the standard manner \cite{Burago}, so that edges in $\cG$ are parametrized in such a way that each of them has length $1$. Associated to $\ell$, there is a complete length metric $d$ on $\cG$ defined by \cite{Burago} 
  \begin{align*} D(x,y) & = \inf_\delta \{~ \ell(\delta) ~ |~\delta:[0,1] \to \cG, ~ \delta(0) = x, ~ \delta(1) = y ~\}. \end{align*}
Thus $(\cG,D)$ becomes a \emph{length} metric space.

\subsection{Space of pointed trees with Gromov-Hausdorff metric}\label{subsec-GH} We call a subgraph $T \subset \cG$ an \emph{infinite tree} if it is non-compact, connected and simply connected. The last condition implies that any loop in $T$ is homotopic with fixed end-points to a trivial loop. Let $X$ be the set of all infinite trees in $\cG$ containing the identity $e \in G$. As a subset of $\cG$, a tree $T\in X$ has an induced length structure $\ell$ and an induced metric $D$. The induced metric $D$ need not be geodesic, since $T$ need not contain shortest paths between its points. However, $D$ induces a length structure $\ell'$ on $T$ which coincides with the restriction of the length structure $\ell$ from $\cG$, and they define a length metric $d$ on $T$. The pair $(T,e) \in X$ with metric $d$ is a pointed metric space, or a pointed tree.

We say that two (possibly finite) pointed graphs $(T,e)$ and $(T',e)$ are\emph{ isomorphic} if and only if there exists a isometry $(T,e) \to (T',e)$ which preserves the distinguished vertex and labeling of edges.

The space $X$ of pointed trees can be given the Gromov-Hausdorff metric $d_{GH}$ \cite{Burago}. One of the metrics with which the Gromov-Hausdorff topology on $X$ can be realised is the ball metric defined now.
 
\begin{defn}\label{defn-ballmetric}\cite{ALM2008,Ghys1999,LozanoSp,Lozano2009}
Let $X$ be the set of all infinite pointed trees in a locally compact Cayley graph $\cG$. Let $T,T' \in X$, and define the distance between $T$ and $T'$ by
  \begin{align} \label{eq-metric1} d_{X} (T,T') & = e^{-r(T,T')}, &   r(T,T') & = \max \{ r \in \mN \cup \{0\} ~ |~ \exists ~ {\rm isomorphism} ~ B_T(e,r) \to  B_{T'}(e,r) \}. \end{align}
The metric $d_X$ is called the \emph{ball metric}.
\end{defn}

\begin{defn}\label{defn-convergence}
Let $\cT =\{(T_n,e)\}_{n=1}^\infty$ be a sequence in $(X,d_X)$. Then $\cT$ converges to $(S,e) \in X$ if and only if for every $r>0$ there exists $n_r >0$ such that for any $n>n_r$ there is an isomorphism $B_{T_n}(e,r) \to B_S(e, r)$. 
\end{defn}

Using the ball metric $d_X$ on $X$ and the fact that $G$ is finitely generated one can prove the following propositions.

\begin{prop}\cite{Ghys1999,LozanoSp}\label{thm-cantorset}
The metric space $(X,d_X)$ is compact and totally disconnected. 
\end{prop} 

If, in addition, $X$ is does not have isolated points, then $X$ is a Cantor set. In particular, if $G = F_n$ or $G = \mathbb{Z}^n$, $n \geq 2$, then $X$ is a Cantor set.

\subsection{Pseudogroup action on the space of pointed trees}\label{subsubsec-actiongraphs} We define a partial action of $F_n$ on the space of pointed trees $X$, which gives rise to a pseudogroup $\fG$ on $X$.

Let $\cP_e(T)$ be the set of paths $\delta:[0,1] \to T$ such that $\delta(0) = e$, $\delta(1) = g \in V(T)$ and $\delta$ is the shortest path between $e$ and $g$ in $T$. The image of $\delta$ in $T$ is the union of edges 
   \begin{align*}w_{h_{i_1}}\cup w_{h_{i_2}}\cup \cdots \cup w_{h_{i_m}} ~ \textrm{where} ~  h_{i_k} \in G_0 ~ \textrm{for} ~ 1 \leq k \leq m.\end{align*}
Thus $\delta$ defines a word $\wth_{i_1} \wth_{i_2} \cdots \wth_{i_m} \in F_n$, where
  \begin{align*}  \left\{\begin{array}{ll}\wth_{i_k} = h_{i_k},  & \textrm{if} ~ \delta^{-1}(s(w_{h_{i_k}})) < \delta^{-1}(t(w_{h_{i_k}})), \\ \wth_{i_k} = h_{i_k}^{-1} , & \textrm{if} ~ \delta^{-1}(s(w_{h_{i_k}})) \geq \delta^{-1}(t(w_{h_{i_k}})). \end{array} \right. \end{align*} 
We note that $g = \wth_{i_1} \wth_{i_2} \cdots \wth_{i_n}$ is a representation in the set $G_0$ of generators of $G$. The procedure we have just described defines an injective map
  \begin{align*} p &: \cP_e(T) \to F_n : \delta \mapsto \wth_{i_1} \wth_{i_2} \cdots \wth_{i_m}, \end{align*}
and the action of $F_n$ on $X$ is defined as follows.

\begin{figure}
\centering
\includegraphics [width=6cm] {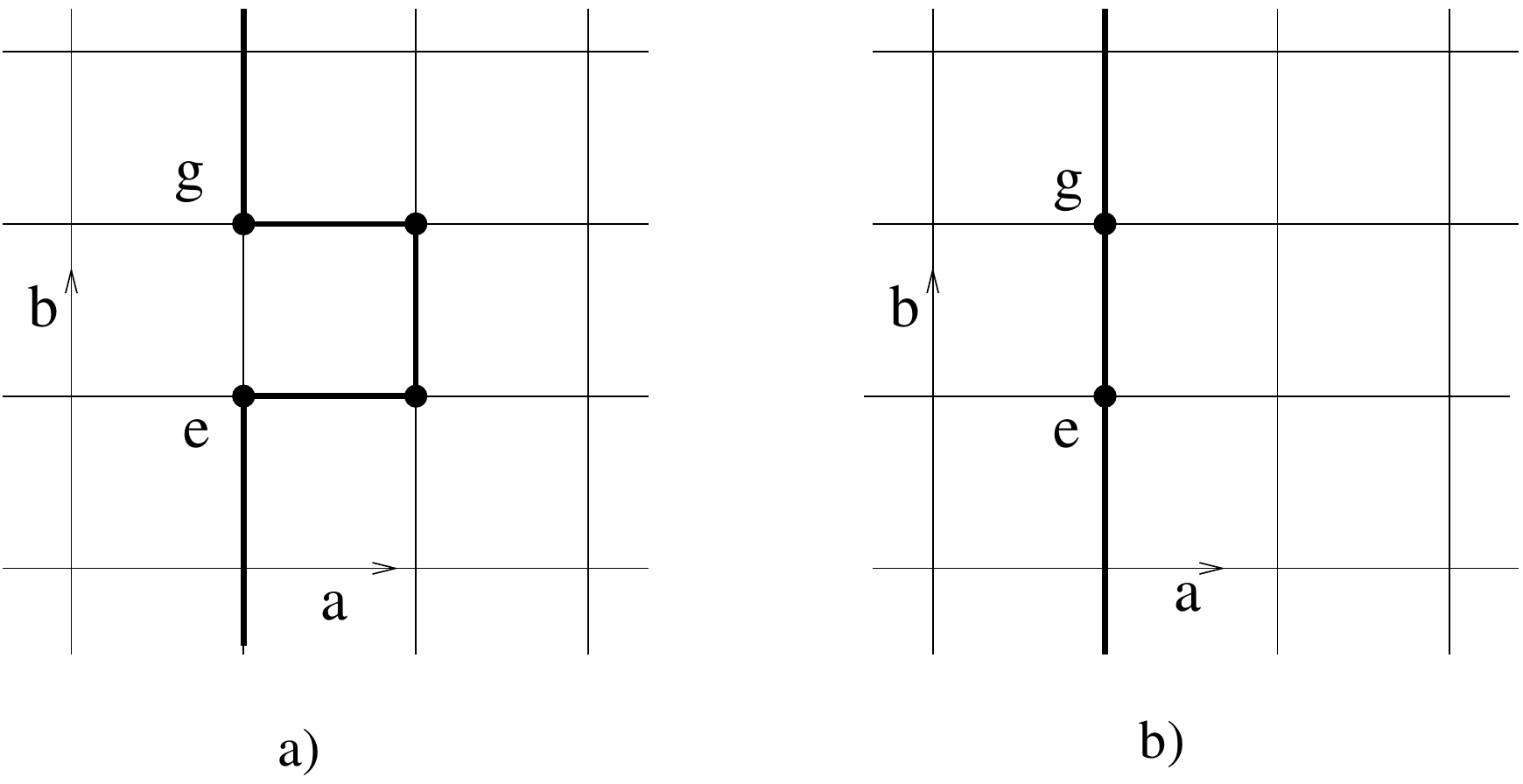}
\caption{The action of $F_2$ on subtrees of the Cayley graph of $\mZ^2$ with a set of generators $G_0 = \{a,b\}$, where $a$-edges are directed to the right, $b$-edges are directed upwards: a) $(T,e) \cdot aba^{-1}$ is defined, $(T,e) \cdot b$ is not defined, b) $(T,e) \cdot b$ is defined, $(T,e) \cdot aba^{-1}$ is not defined. }
 \label{fig:Ps-def}
\end{figure}

\begin{defn}\label{defn-actionany}
Let $n < \infty$ be the cardinality of a set $G_0$ of generators of $G$, and $(X,d_X)$ be the corresponding set of pointed trees. An action of $g \in F_n$ on $X$ is defined as follows.
\begin{enumerate}
\item $(T,e) \cdot g$ is defined if and only if there exists a path $\delta \in \cP_e(T)$ such that $p(\delta) = g$.
\item $(T',e) = (T,e) \cdot g$ if and only if there is an isomorphism of pointed spaces $\alpha:(T,g) \to (T',e)$. 
\end{enumerate}
\end{defn}

To a partial action of $F_n$ on $X$ we can associate a pseudogroup of local homeomorphisms $\fG$ as follows. For $r>0$ denote by $D_{X}(T,r)$ a clopen subset of diameter $e^{-r}$ about $(T,e)$, that is,
  \begin{align*}D_{X}(T,r) & = \left\{ (T',e) ~ | ~ d_{X}(T,T') \leq e^{-r} \right\}. \end{align*}
For each $g \in F_n$ let $\ell_g = d_{\cF_n}(e,g)$. The action of $g$ is defined on the union of clopen subsets
  \begin{align}\label{domgamma}D & = \bigcup \left\{D_X(T,\ell_g) ~ |~ T \in X, ~ \delta \in \cP_e(T) ~ \textrm{such that} ~ p(\delta) = g \right\}, \end{align}
which is clopen since $G_0$ is a finite set and so \eqref{domgamma} is a finite union. The mapping
  \begin{align*} \gamma_g & : {\rm dom}(\gamma_g) = D \to X_n : (T,e) \mapsto (T,e)\cdot g \end{align*}
is a homeomorphism onto its image, and a pseudogroup $\fG_n$ is defined to be a collection  
  \begin{align}\label{eq-pseudogroup}\fG_n& = \bigsqcup \{ ~\gamma_g ~|~ g \in F_n ~ \}\end{align}
of local homeomorphisms. A subset $\fG^0_n = \{~\gamma_g \in \fG ~|~ g \in G_0 ~\}$ is a generating set of $\fG$.  

\begin{defn}\label{defn-eqrel}\cite{Blanc2001}
Let $(X,d_X)$ be a space of pointed trees, and $\fG$ be the pseudogroup on $X$. Then $(T,e)$ and $(T',e)$ are \emph{$\cR$-equivalent} 
  $$(T,e) \sim_{\cR} (T',e)$$
if and only if there exists $\gamma_g \in \fG$ such that $\gamma_g(T,e) = (T',e)$. An equivalence class of $(T,e) \in X$ with respect to $\cR$ is denoted by $\cR(T)$ and is called the \emph{orbit} of $(T,e)$ under the action of $\fG$.
\end{defn}

\subsection{Foliated space $\fM_G$ with foliation by Riemann surfaces}\label{subsec-folRS} We realize the pseudogroup dynamical system $(X,\fG)$ as the holonomy system of a smooth foliated space $\fM_G$.

\begin{thm}\cite{Ghys1999,LozanoSp}\label{thm-genconstruction}
Let $G$ be a finitely generated group, and $(X,d_X)$ be the corresponding space of pointed trees with the action of a pseudogroup $\fG$. Then there exists a compact metric space $\fM_G$, and a finite smooth foliated atlas $\cV = \{\phi_i:V_i \to U_i \times \fX_i\}_{1 \leq i \leq \nu}$, where $U_i \subset \mR^2$ is open, with associated holonomy pseudogroup $\fP$, such that the following holds. 
\begin{enumerate}
\item The leaves of $\fM_G$ are Riemann surfaces.
\item There is a homeomorphism onto its image
  \begin{align*} t: X \to \cup_{1 \leq i \leq \nu} \fX_i, \end{align*}
such that $t(X)$ is a complete transversal for the foliation $\cF$, and $\fP|_{\tau(X)} = t_*\fG$, where $t_*\fG$ is the pseudogroup induced on $t(X)$ by $\fG$.
\end{enumerate}
\end{thm}

For completeness we give a sketch of the proof \cite{Ghys1999,LozanoSp}.

\emph{Sketch of proof}. Let ${\bf A}$ be the finite set of connected subtrees of an open ball $B_{\cG}(e,1)$, excluding the subgraph consisting of a single point. Denote by $n_{\bf a}$ the number of edges of ${\bf a} \in {\bf A}$. Let $\Sigma_{\bf a}$ be a compact surface with boundary homeomorphic to a $2$-sphere with $n_{\bf a}$ disks taken out, and label the connected components of the boundary as follows: each boundary component corresponds to a labeled edge $w_h$ of ${\bf a}$, a boundary component is labeled by $h$ if $s(w_h) =e$, and it is labeled by $h^{-1}$ if $t(w_h) = e$. Choose a Riemannian metric on $\Sigma_{\bf a}$ in such a way that each connected component of the boundary has a closed neighborhood $Y_g^{\bf a}$ isometric to $\mR/\mZ \times [0,1/2]$, and fix these metrics. For each ${\bf a} \in {\bf A}$ choose a base point $p_{\bf a}$ in the interior of $\Sigma_{\bf a}$.
 
Form a disjoint union $\bigsqcup_{{\bf a} \in {\bf A}} D_X({\bf a},1) \times \Sigma_{\bf a}$, and identify neighborhoods of boundary components as follows. Suppose $(T,e) \in D_X({\bf a},1)$ and $\Sigma_{\bf a}$ contains a boundary component marked by $h \in G_0 \cup G_0^{-1}$, where $G_0^{-1} = \{h^{-1} ~|~ h \in G_0\}$. Then there is a pointed tree
 $$(T',e) = (T,e) \cdot h \in D_X({\bf a}',1),$$ 
and $\Sigma_{\bf a'}$ has a boundary component marked by $h^{-1}$. We identify $\{(T,e)\} \times Y^{\bf a}_h \subset \{(T,e)\} \times \Sigma_{\bf a}$ and $\{(T',e)\} \times Y^{\bf a'}_{h^{-1}}\subset \{(T',e)\} \times \Sigma_{\bf a'}$ by setting
   $$\left((T,e),(\theta,s)\right) \sim \left((T',e),(-\theta, 1/2 - s)\right),~ \theta \in [0,1], ~ s \in [0,1/2].$$
Taking the quotient by this equivalence relation one obtains a space
  \begin{align*}\fM_G &: =\bigsqcup_{{\bf a} \in {\bf A}}  D_X({\bf a},1) \times \Sigma_{\bf a}/\sim. \end{align*}
The obtained space $\fM_G$ is compact foliated space \cite{Ghys1999, LozanoSp}, with leaves Riemann surfaces. Given a clopen cover $\cZ_{\bf a}$ of $D_X({\bf a},1)$, and a geodesically convex open cover $\cU_{\bf a}$ of a surface $\Sigma_{\bf a}$, by taking products of charts one constructs the required foliation cover $\cV$, with transverse space given by
\begin{align*} \fX & = \bigsqcup \{ Z^{\bf a}_j ~|~ j \in J_{\bf a}, ~ {\bf a} \in {\bf A}~\}. \end{align*}
To obtain an embedding of $X$ notice that there is an embedding
  \begin{align} \label{eq-emb}\imath &: X \to \fM_G : (T,e) \mapsto \imath_{\bf a} \left((T,e), p_{\bf a}\right), \end{align}
and define $t:X \to \fX$ as the obvious composition of \eqref{eq-emb} with the chart maps of $\cV$. Define a metric $d_{\fX}$ on $\fX$ as follows: since each $Z^{\bf a}_j$ is a subset of $X$ we can define the metric $d_{\fX}$ on $Z^{\bf a}_j$ by restricting the metric $d_X$. If $w \in  Z^{\bf a}_j$, $w' \in  Z^{\bf a'}_{j'}$ and $w,w' \in t(X)$, then set 
  $$d_{\fX}(w,w') = d_X(t^{-1}(w),t^{-1}(w')),$$ 
otherwise set $d_{\fX}(w,w') =1$. The properties of $t$ stated in the formulation of the theorem follow straightforwardly from the construction.
\endproof

\subsection{Universal property} Theorem \ref{prop-univprop} is a consequence of Theorem \ref{thm-genconstruction}. Given a pointed tree $(T,e) \in X$, we denote by $L_T$ a leaf such that $ \imath(T,e) \in L_T$.

\proof \emph{(of Theorem \ref{prop-univprop})}. Denote by $n$ the cardinality of the generating set $G_0$ of $G$. Identify $G_0$ with a generating set of a free group on $n$ generators $F_n$.  Let $(X,d_{X})$ be a space of infinite pointed trees contained in the Cayley graph $\cG$ of $G$, and $(X_n,d_{X})$ be a space of infinite pointed trees contained in the Cayley graph $\cF_n$ of $F_n$. Then every tree in $\cG$ can be identified with a unique tree in $\cG_n$, and so $X$ is identified with a subset of $X_n$. Thus there is an embedding $\Phi:X \to X_n$ such that $\Phi(X)$ is closed in $X_n$. It is then clear that the restricted pseudogroup $\Phi_*\fG = \fG_n|_{\Phi(X)}$, and the result follows as a consequence of Theorem \ref{thm-genconstruction}.
\endproof

\subsection{Graph matchbox manifolds} \label{subsec-recurrence}

We give a definition of a graph matchbox manifold.

\begin{defn}
Let $G$ be a finitely generated group, and $\fM_G$ be a smooth foliated space obtained by the construction of Kenyon and Ghys. Then a \emph{graph matchbox manifold} is the closure $\cM = \overline{L}$ of a leaf $L$ in $\fM_G$.
\end{defn}

Given a pointed tree $(T,e) \in X$, we denote by $L_T$ a leaf such that $\imath(T,e) \in L_T$. It follows from the proof of Theorem \ref{thm-genconstruction} that $L_T \cap \imath(X) = \imath(\cR(T))$, where $\cR(T)$ denotes the orbit of $(T,e) \in X$ under the action of the pseudogroup $\fG$. As usual in foliation theory, one aims to relate asymptotic properties of orbits in $X$ with asymptotic properties of leaves in $\fM_G$, and thus study the dynamics of graph matchbox manifolds via orbits of points in $X$. Using the construction in \cite[section 2.2]{ClarkHurder2010b} of a cover well-adapted to the metrics $d_\fM$ and $d_\fX$, and uniform metric estimates which measure distortions between $d_\fM$ and $d_\fX$ \cite[section 2.2]{ClarkHurder2010b}, one obtains the following lemma.

\begin{lemma}\label{defn-graphMM}
$\cR(T') \subset \overline{\cR(T)}$ if and only if $L_{T'} \subset \overline{L_T}$.
\end{lemma}

In particular, if $(T,e) \subset X$ is a pointed tree in $X$, and $\cM_T = \overline{L_T}$ is the corresponding graph matchbox manifold, then
  \begin{align*} \cM_T & = \left\{ L \subset \fM_G ~ | ~ \imath\left(\overline{\cR(T)} \right) \cap L \ne \emptyset \right\}.\end{align*}

As a consequence of Lemma \ref{defn-graphMM} one obtains the relation between asymptotic properties of leaves in $\cM_T$ and the topology of the set $\overline{\cR(T)}$. We first recall some definitions.

\begin{defn}\label{limsetleaf}\cite{CandelConlon2000}
Let $\{K_\alpha\}$ be the set of all compact subsets of $L_{T}$ and $W_\alpha = L_{T} - K_\alpha$ be the complement of $K_\alpha$ in $L_{T}$. Then the \emph{limit set}, or the \emph{asymptote}, of $L_{T}$ is the set 
  \begin{align*}\lim L_{T} & = \bigcap_{\alpha} \overline{W_\alpha}. \end{align*}
A leaf $L_T$ is \emph{recurrent} if $L_T \subset \lim L_T$. A graph matchbox manifold $\overline{L_T}$ is \emph{recurrent} if $L_T$ is recurrent.
\end{defn}

In a compact foliated space $\fM_G$ the limit set $\lim L$ is compact, non-empty and saturated \cite{CandelConlon2000}. The following corollary follows straightforwardly from Lemma \ref{defn-graphMM} and from definitions.

\begin{cor}\label{recurrencelemma}
A graph matchbox manifold $\cM_T = \overline{L_T}$ is recurrent if and only if $\overline{\cR(T)}$ is a Cantor set. If a matchbox manifold $\cM$ has a recurrent leaf, then it contains an uncountable number of leaves. If a transitive leaf $L \subset \fM$ is not recurrent, then every point of $\cR(T)$ is isolated in $\overline{\cR(T)}$.
\end{cor}

Finally it is convenient to formulate the following criterion of recurrence in terms of balls in graphs.

\begin{lemma}\label{lemma-recurrencecritgraph}
A graph matchbox manifold $\cM_T$ is recurrent if and only if there exists a sequence $\{g_\ell\} \in V(T)$, $g_\ell \ne e$, such that for every $r>0$ there exists $\ell_r>0$ such that for all $\ell \geq \ell_r$ there is an isomorphism 
  $$\alpha^r_\ell: D_T(e,r) \to D_T(g_\ell,r).$$
\end{lemma}

\subsection{Ends of leaves}\label{section-endsofgraphs}

In this technical section, we make precise the relationship between the ends of a graph $T$ and the ends of the leaf $L_T$. A tool for that is the \emph{graph of the $\fG$-orbit} of $(T,e) \in X$ which was considered in \cite{LozanoSp}, and which we define now.

\begin{defn}\label{defn-holgraph}
Let $(T,e) \in X$. A graph $\G_T$ of the $\fG$-orbit $\cR(T)$ is defined as follows: we set $V(\G_T) = \cR(T)$, and we join $(T',e)$ and $(T'',e) \in \cR(T)$ by an oriented edge $w$ with $s(w) = (T',e)$ and $t(w) = (T'',e)$ if and only if there exists $\gamma_h \in \fG^0$ such that $\gamma_h(T',e) = (T'',e)$. 
\end{defn}

The graph $\G_T$ is given a length structure in the usual manner so that the length of each edge is $1$. 

Ends are the means to study asymptotic properties of a topological space $S$, and are a form of compactification of $S$. The following definition is also convenient for computations.

\begin{defn}\cite{CandelConlon2000}
Let $S$ be a Hausdorff locally compact locally connected separable topological space. Let $K_1 \subset K_2 \subset \ldots$ be an increasing sequence of compact subsets $K_i$ such that $\bigcup_i K_i = S$, and 
  $$\{U_{\alpha_i}\} : U_{\alpha_1} \supset U_{\alpha_2} \supset \ldots,$$ 
be a decreasing sequence, where $U_{\alpha_i}$ is an unbounded connected component of $S \setminus K_i$. We say that $\{U_{\alpha_i}\}$ and $\{V_{\beta_i}\}$ are equivalent, $\{U_{\alpha_i}\} \sim \{V_{\beta_i}\}$, if for every $i>0$ there exists $j>i$ such that
  $$\left( U_{\alpha_i} \cap V_{\beta_i} \right) \supset \left( U_{\alpha_j} \cup V_{\beta_j} \right).$$
An \emph{end} ${\bf e}$ of $S$ is an equivalence class of $\{U_{\alpha_i}\}$ with respect to the equivalence relation $\sim$.
\end{defn}

We denote by $\cE(S)$ the set of ends of $S$ and topologize $S^* = S \cup \cE(S)$ as follows \cite{CandelConlon2000}. For every ${\bf e}$ and $U_{\alpha_i} \in \{U_{\alpha_i} \}$, an open set $U_{\alpha_i}$ is called the fundamental neighborhood of ${\bf e}$. The topology on $S^*$ is given by all open sets of $S$, plus all fundamental neighborhoods of ends together with ends contained in each neighborhood. We say that a sequence $\{x_n\} \in S^*$ converges to an end ${\bf e} \in \cE(S)$ if every fundamental neighborhood of ${\bf e}$ contains all but a finite number of elements of $\{x_n\}$.

It was proved in \cite[Section 3.2]{CandelConlon2003}, that the number of ends of a leaf in a foliated space is the same as the number of ends of a corresponding path-connected component of the holonomy graph of the foliated space. A very similar argument, relying on the existence of embeddings 
  $$\imath_{\bf a}:D_X({\bf a},1) \times \Sigma_{\bf a} \to \fM_G, ~ {\bf a} \in {\bf A},$$ 
which provide a cover of $\fM_G$, and which allow us to associate to each orbit of $\fG$ in $X$ a `plaque chain' in $L_T$ made up of compact subsets isometric to $\Sigma_{\bf a}$, ${\bf a} \in {\bf A}$, one obtains the following conclusion.

\begin{lemma}\label{lemma-holgraph}
Let $L_T \in \fM_G$ be a leaf, and $\G_T$ be the graph of $\cR(T)$. Then there is a homeomorphism $\cE(L_T) \to \cE(\G_T)$.
\end{lemma}

Define a map $v_T : T \to \G_T$  by $v_T(g) = (T,e) \cdot g$ on the sets of vertices $V(T)$ and $V(\G_T)$, and map an edge $m \in E(T)$ isometrically onto the oriented edge starting at $v_T(s(m))$ and ending at $v_T(t(m))$. 

\begin{lemma}\label{lemma-projection}
The map $v_T:T \to \G_T$ is a covering projection. If $v_T$ is a homeomorphism, there is an induced homeomorphism $v_T^*:\cE(T) \to \cE(\G_T)$ of end spaces.
\end{lemma}

\subsection{Examples} \label{examples-prel} Let $G = F_2$ and $G_0 = \{a,b\}$. As before, the Cayley graph of $F_2$ is denoted by $\cF_2$.

\emph{Genus two surface and tori}. The graph $(\cF_2,e)$ is invariant under the action of the pseudogroup $\fG_2$, that is, for any $g \in G_2$ we have $\gamma_g(\cF_2,e) = (\cF_2,e)$. It follows that $L_{\cF_2}$ is homeomorphic to a genus two surface. 

\begin{figure}
\centering
\includegraphics [width=3.6cm] {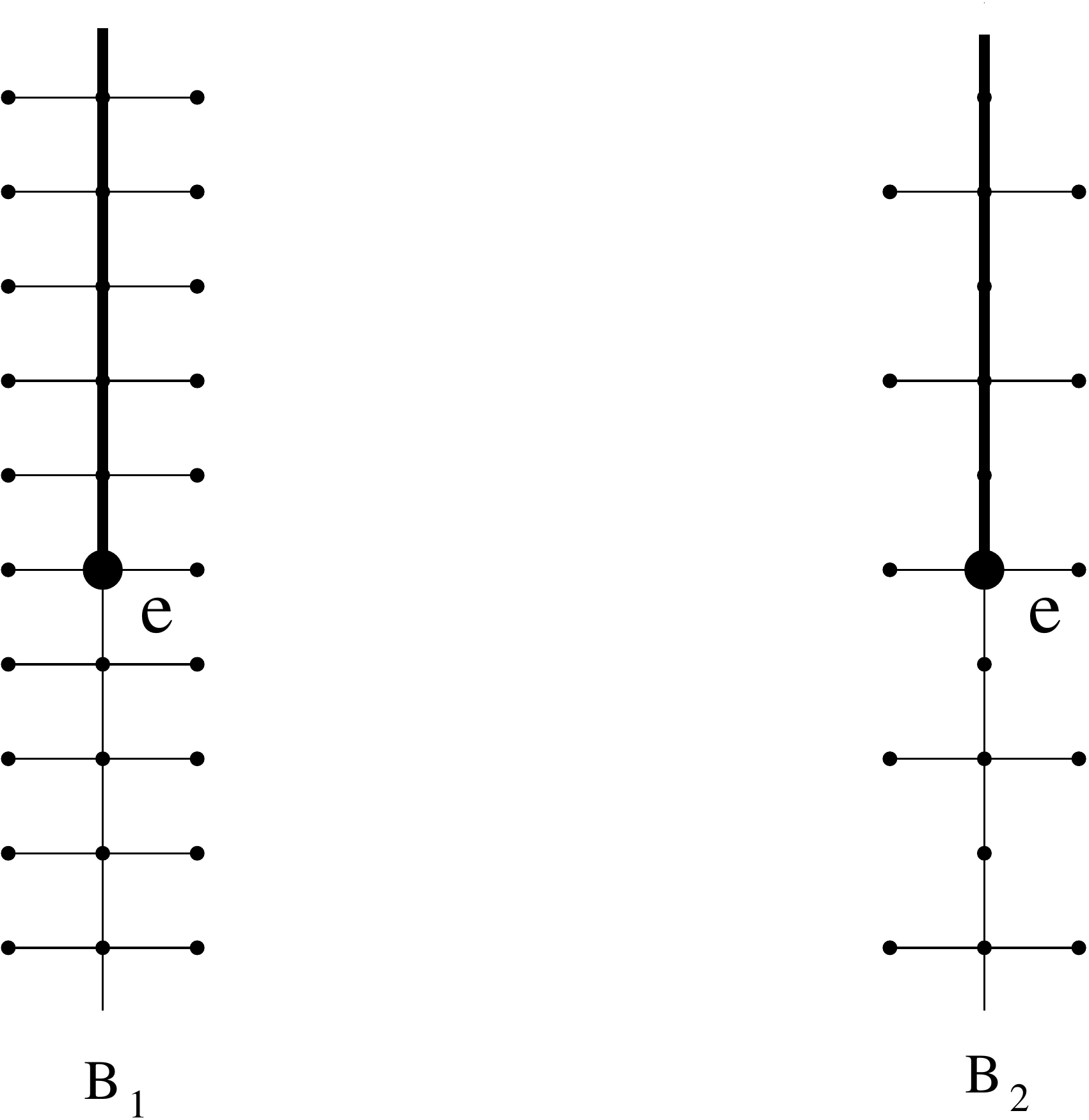}
\caption{Graphs $B_1$ and $B_2$. Thick lines indicate paths $\delta_1$ and $\delta_2$ in $B_1$ and $B_2$ respectively.}
 \label{fig:graphsBn-1}
\end{figure}

Let $B_0 \subset \cF_2$ be a subgraph with the set of vertices $V(B_0) = \bigcup_{n \geq 0} \{b^n,b^{-n}\}$. Then $(B_0,e) \in {\rm dom} (\gamma_{b^{\pm n}})$, and $\gamma_{b^{\pm n}}(B_0,e) = (B_0,e)$, and $L_{B_0}$ is homeomorphic to the standard torus. Leaves $L_{B_1}$ and $L_{B_2}$ corresponding to bi-infinite graphs $B_1$ and $B_2$ in Fig. \ref{fig:graphsBn-1} are also homeomorphic but not isometric to the standard torus.

\emph{Fusion of two graphs}. We now create a graph $T$ such that the corresponding leaf $L_T$ accumulates on $L_{B_1}$ and $L_{B_2}$. The graph $T$ is obtained as a \emph{fusion} of $B_1$ and $B_2$. Here we only describe the idea of the fusion, with the rigorous algorithm given in Section \ref{section-directed}. Fusion plays an important role in the proof of main theorems in Sections \ref{asymptotic} and \ref{exmpl}.

\begin{figure}
\centering
\includegraphics [width=10.5cm] {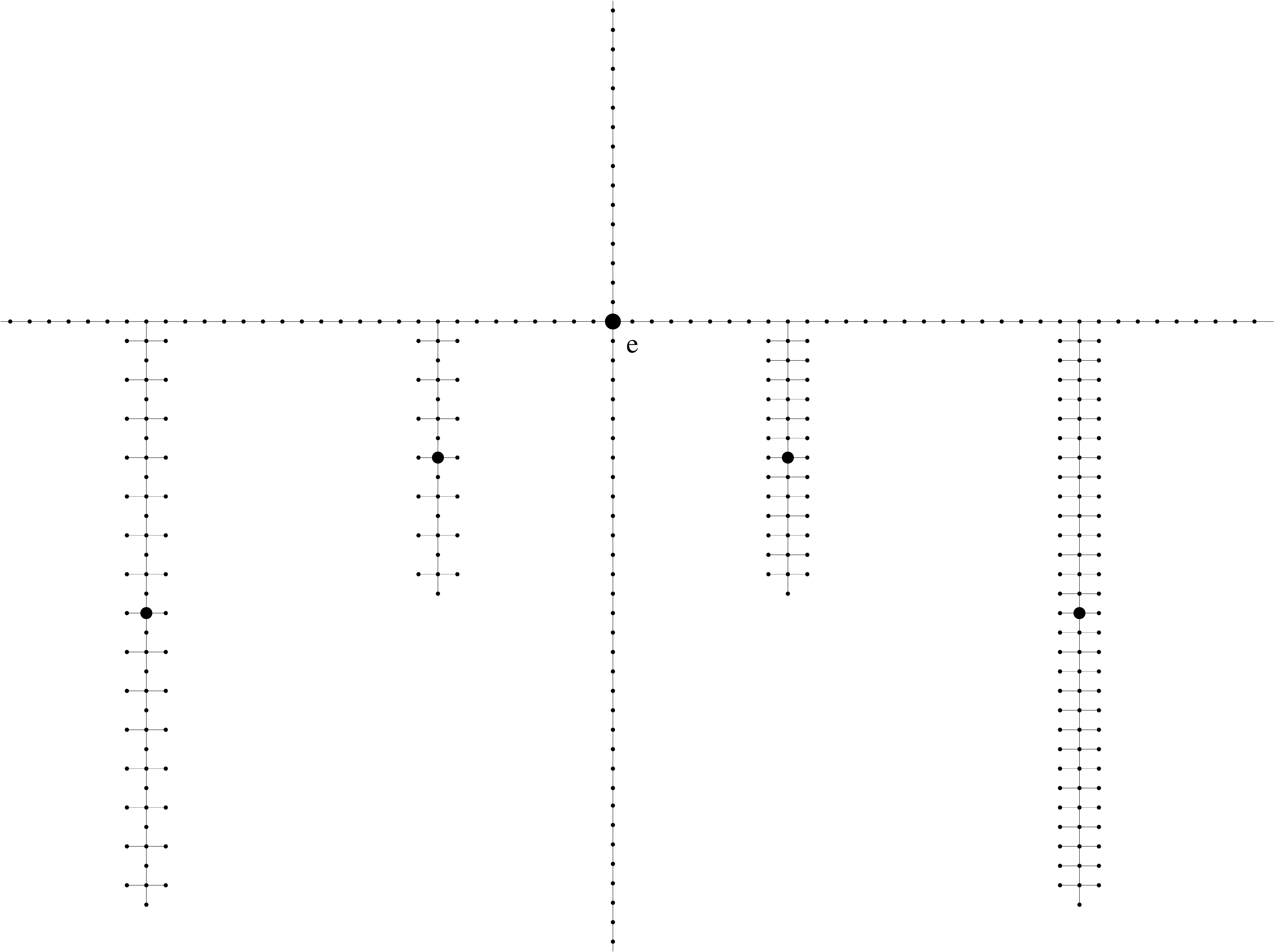}
\caption{Fusion of graphs $B_1$ and $B_2$: infinite lines of $a$- and $b$-edges intersecting at the origin $e$ with isomorphic copies of $D_{B_k}(e,7)$ and $D_{B_k}(e,15)$, $k=1,2$, attached. Large unmarked dots indicate images of the centers $e$ of these balls under the attaching maps.}
 \label{fig:fusion}
\end{figure}

We start with a graph $T' \subset \cF_2$ such that $V(T') = \bigcup_{n \geq 0} \{a^n, a^{-n},b^n,b^{-n}\}$ ($T'$ is a union of infinite lines of $a$-edges and of $b$-edges, intersecting at $e$). The idea is to attach to $T'$ copies of $D_{B_1}(e,r)$ and $D_{B_2}(e,r)$, $r>0$, in a specified order to obtain a graph $T$, so that the graph $\G_T$ of the orbit $\cR(T)$ has $4$ ends, there is at least one end accumulating on $(B_1,e)$, and at least one end accumulating on $(B_2,e)$.

\begin{figure}
\centering
\includegraphics [width=10.5cm] {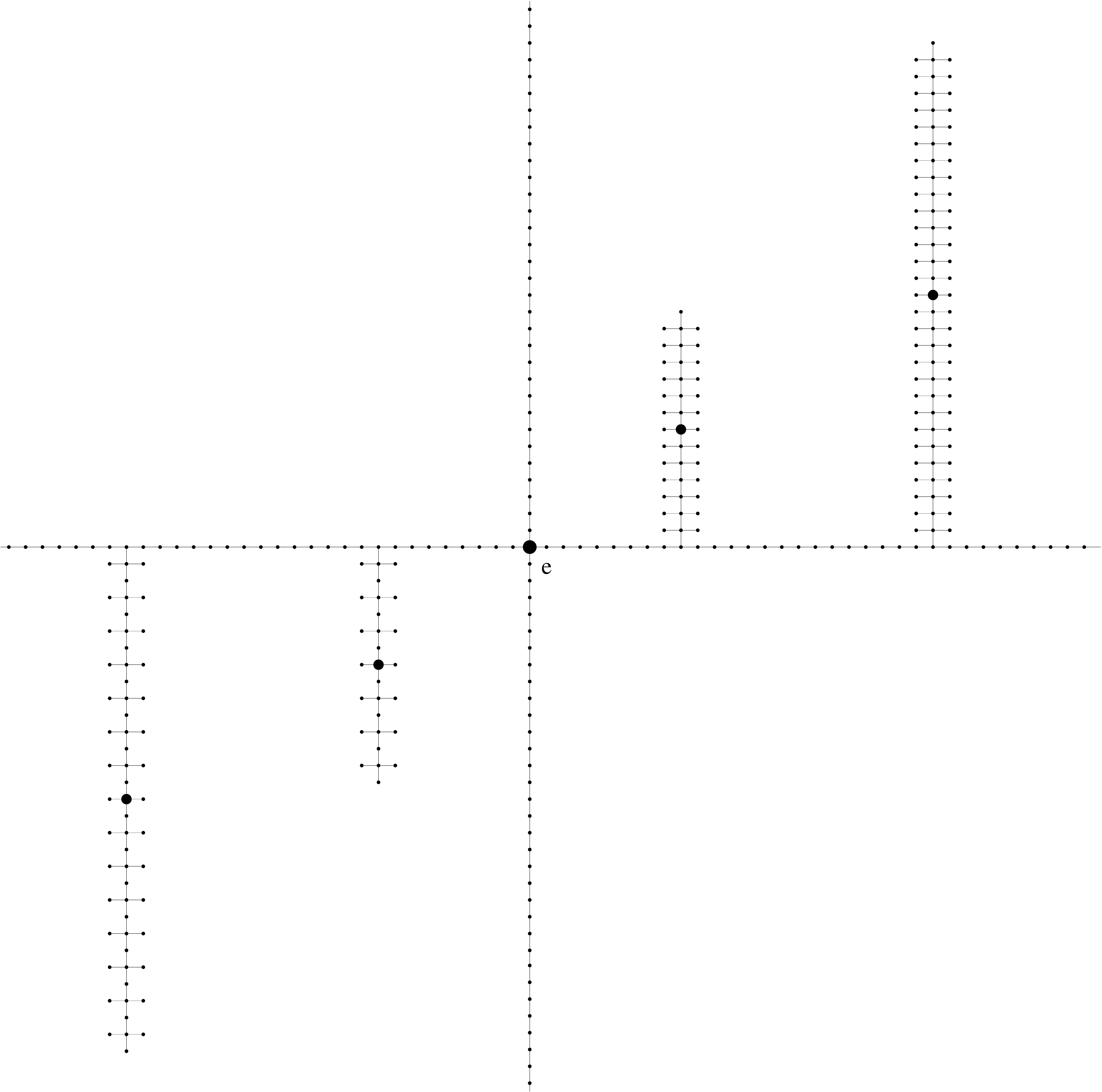}
\caption{Fusion of graphs $B_1$ and $B_2$, for a different choice of paths $\delta_k$ in $B_k$, $k=1,2$. Large unmarked dots indicate images of the centers $e$ of these balls under the attaching maps.}
 \label{fig:fusion-2}
\end{figure}

Let $r_0 = 2$, and notice that $D_{B_1}(e,2) \ne D_{B_2}(e,2)$. Let $\delta_1$ and $\delta_2$ be infinite subgraphs in $B_1$ and $B_2$ respectively, such that $V(\delta_i) = \bigcup_{n \geq 0}\{b^n\}$ (see Fig. \ref{fig:graphsBn-1}). 

For $i>0$ let 
  $$r_i = 2^{r_0+i} -1 ~ \textrm{and} ~ R_i = \sum_{k=0}^{i-1} r_k.$$ 
For $k=1,2$ let $w_i^k$ be an edge of $\delta_k$ such that $w_i^k \subset D_{B_k}(e,r_i)$ and $w_i^k \cap \partial D_{B_k}(e,r_i) \ne \emptyset$. Then $T$ is obtained by attaching to $T'$ isomorphic copies of $D_{B_1}(e,r_i)$, $i \in \mathbb{N}$, such that the image of $w_i^1$ has a common vertex with the axis of $a$-edges, and this vertex is $a^{R_i +r_i}$; and also attaching to $T'$ isomorphic copies of $D_{B_2}(e,r_i)$, $i \in \mathbb{N}$, such that the image of $w_i^2$ has a common vertex with the axis of $a$-edges, and this vertex is $a^{-(R_i +r_i)}$. The orbit $\cR(T)$ of the resulting graph $(T,e)$ accumulates on $(B_1,e)$ and $(B_2,e)$. The closure $\overline{R(T)}$ also contains other pointed trees, for example, $(B_0,e) \subset \overline{R(T)}$. It follows that $\overline{L}_T \supset \{L_{B_0},L_{B_1},L_{B_2}\}$.

We note that the resulting graph $T$ depends on the choice of infinite paths $\delta_k$ in $B_k$, $k=1,2$. For example, if $V(\delta_1) = \bigcup_{n \leq 0}\{b^n\}$ in $(B_1,e)$, then by fusion we obtain the graph in Fig. \ref{fig:fusion-2}.

\begin{figure}
\centering
\includegraphics [width=9.0cm] {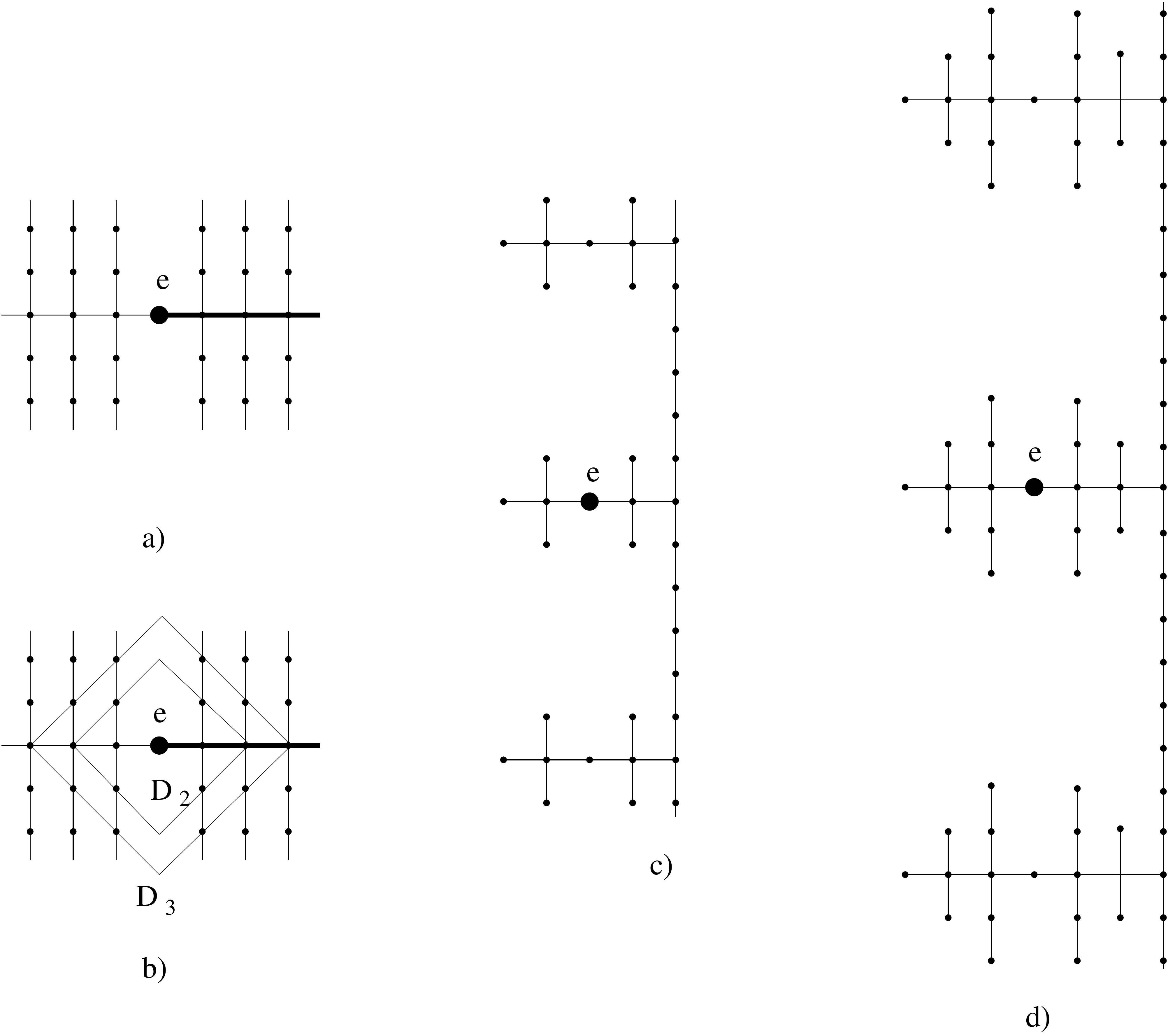}
\caption{Constructing graphs with periodic orbits: a) Graph $T$.  Thick line indicates a path $\delta$. b) Closed balls $D_2$ and $D_3$ in the graph $T$. c) Graph $S_2$. d) Graph $S_3$.}
 \label{fig:compact}
\end{figure}

\emph{Graphs with periodic orbits in a neighborhood of a given graph}. Let $\delta_r = e^{-r}$ and let $T$ be a graph with the set of vertices (see Fig. \ref{fig:compact}, a))
  \begin{align*} V(T) & = \{e\}  \bigcup_{i \geq 1} \left( \bigcup_{j \geq 0} \{a^ib^j,~a^ib^{-j},a^{-i}b^j,a^{-i}b^{-j}\} \right).  \end{align*}
The idea is to attach copies of $D_r = D_T(e,r)$ to a two-ended subgraph of $\cF_2$ at regular intervals. 

 For each $r>0$ there is an edge $w_r$ labeled by $a$, such that $w_r \subset D_T(e,r)$ and $w_r \cap \partial D_T(e,r) \ne \emptyset$. For instance, we can choose $w_r$ as the edge lying in the intersection of $D_T(e,r)$ and an infinite edge path $\delta$ in $(T,e)$ such that $V(\delta) = \bigcup_{n \geq 0} \{a^n\}$ (see Fig. \ref{fig:compact}, a)). Let $S_r'$ be a subgraph of $\cF_2$ with vertices $V(S_r') = \bigcup_{i \geq 0} \{a^rb^i,~ a^rb^{-i}\}$. For each $i \in \mathbb{Z}$, attach to $S_r'$ an isomorphic copy of $D_r$ (see Fig. \ref{fig:compact} b)) such that the edge $w_r$ has a common vertex $a^rb^{3i}$ with $S'$. The resulting graph $S_r$ is periodic, and so the leaf $L_{S_r}$ is compact. We note that $d_{X_n}(T,S_r) = e^{-r}$ (see Fig. \ref{fig:compact} b), c)).

\section{Hierarchy of graph matchbox manifolds}\label{asymptotic}

Let $\fM_G$ be a foliated space obtained by the construction of Kenyon and Ghys, and $\cS$ be the set of graph matchbox manifolds, that is, the set of transitive saturated closed subsets of $\fM_G$. In this section we study a partial order $\preceq$ on $\cS$ given by inclusions, and prove Theorem \ref{thm-globalleaf}.

\begin{defn}\label{defn-partialorder}
Let $\fM_G$ be a foliated space and $\cS$ be the set of transitive saturated closed subsets in $\fM_G$. A partial order $\preceq$ on $\cS$ is given by the following rule:  
  \begin{align*}\textrm{Let} ~ \cM_1,\cM_2 \in \cS, ~\textrm{then} ~ \cM_1 \preceq \cM_2 ~ \textrm{ if and only if}~ \cM_1 \subseteq \cM_2.\end{align*}
\end{defn}

By Theorem \ref{prop-univprop} a partially ordered set $(\cS,\preceq)$ is a subset of $(\cS_n,\preceq)$, where $\cS_n$ is the partially ordered set of graph matchbox manifolds in $\fM_n$, $n = {\rm card}(G_0)$. We can then restrict our study to $\fM_n$. 

\subsection{Compact leaves} We prove Theorem \ref{thm-globalleaf}.(1), that is, the set 
  $$C = \{L \subset \fM_n ~|~ L ~ \textrm{is compact}\}$$ 
is a dense meager subset and $C \cap \cX$ is countable, where $\cX = \imath(X)$ is an embedding of $X$ into $\fM_n$ given by formula \eqref{eq-emb}. The proof is contained in Lemmas \ref{prop-compactdense} and \ref{prop-meager}.

We denote by $\partial U$ the boundary of a subset $U \subset \cF_n$. We notice that if $r>0$ is an integer then the boundary $\partial D_{\cF_n}(e,r)$ contains only vertices.

\begin{lemma}\label{prop-compactdense}
The subset $C = \{L \subset \fM_n ~ | ~ L ~ \textrm{is compact}\}$, $n>1$, is dense in $\fM_n$.
\end{lemma}

\proof By Lemma \ref{defn-graphMM} it is enough to prove that the set $\cC = \{(S,e) \in X ~|~ L_S \in C\}$ is dense in $X$, that is, for every $(T,e) \in X$ and every $\delta>0$ there is an $(S,e) \in \cC$ such that $d_X(T,S) < \delta$. 

Given $\delta> 0 $, choose $r > 0$ such that $e^{-r} < \delta$. The procedure of constructing $S$ consists of attaching copies of a ball $D_T(e,r)$ to a two-ended subgraph of $\cF_n$ at regular intervals (see Section \ref{examples-prel} for an example).

More precisely, since $T$ is infinite, it has at least $1$ end. Then there exists an edge $w_h$, $h \in G_0$, such that $w_h \subset D_T(e,r)$ and $w_h \cap \partial D_T(e,r) \ne \emptyset$. Set
  $$\beta = \left\{\begin{array}{ll} \phantom{-}1, &  t(w_h) \in \partial D_T(e,r), \\ -1, & s(w_h) \in \partial D_T(e,r).\end{array} \right.$$

Let $v_2 = \partial D_T(e,r) \cap w_h$, and let $v_1$ be the other vertex of $w_h$. 

Choose $\widehat{h} \in G_0$ such that $\widehat{h} \ne h$. Let $c:[0,1] \to T$ be the shortest path joining $e \in T$ with $v_2$. Then $w_h \subset c([0,1])$. We associate to $c$ a finite word $\wth_1 \cdots \wth_{k-1}h^\beta \in F_n$, where $k= d_T(e,v_2)$, and either $\wth_s = h_s$ or $\wth_s = h_s^{-1}$ for some $h_s \in G_0$, and $V(c([0,1])) = \{v_2\} \cup \left( \cup_{1\leq s \leq k-1} \wth_s \right)$. Now let $f =  \wth_1 \cdots \wth_{k-1}h^{\beta}$ and
  \begin{align*}g^{n} & = \left(f \widehat{h}^{3r} (f)^{-1}\right)^n,  && g^{-n} = \left(f \widehat{h}^{-3r} (f)^{-1}\right)^n. \end{align*} 
Denote by $\delta^+$ and $\delta^-$ respectively the subgraphs of $\cF_n$ with the vertex sets
  \begin{align*} V(\delta^+) & = \bigcup_{i = 1}^{3r} f\widehat{h}^i, && V(\delta^-) = \bigcup_{i = 1}^{3r} f\widehat{h}^{-i}\end{align*}
that is, $\delta^\pm$ are paths comprised of $3r$ edges marked by $\widehat{h}$ and starting at $f$. The path $\delta^+$ traverses oriented $\widehat{h}$-edges in the positive direction, and $\delta^-$ in the negative direction. Denote by $V_0^\pm = D_T(e,r) \cup \delta_\pm$, and let $V^{\pm n}$ be a subgraph of $\cF_n$ such that there is an isomorphism
   $$\alpha_{\pm n}: (V^{\pm n},g^{\pm n}) \to (V_0^\pm,e). $$ 
Then define
  \begin{align*} S & = (V_0^+ \cup V_0^-) \bigcup \left( \bigcup_{n \in \mZ \setminus\{0\}} V^{n} \right). \end{align*}
The resulting infinite graph is connected and has two ends. By construction $(S,e)$ is invariant under the action of $g^{\pm n}$, $n \in \mN$, and $\G_{S}$ is a finite graph. As desired, $d_X\left((T,e),(S,e) \right) < \delta$. 
\endproof

Recall \cite{Kur} that a set $A \subset M$ is \emph{nowhere dense} if $M \setminus \overline{A}$ is dense in $M$. A subset $A \subset M$ is \emph{meager} (or \emph{of the first category}) if it is the union of a countable sequence of nowhere dense sets.

Let $(T,e) \in X$ be a graph and let $\G_T$ be the graph of $\cR(T)$. We construct a section of the covering projection $v_T:T \to \G_T$ (see Section \ref{section-endsofgraphs}) on the set of vertices $V(T)$ as follows. For every $(T',e) \in V(\G_T)$ let $g_{T'} \in V(T)$ be a vertex such that 
  \begin{align*} d(e,g_{T'}) & = \min \{d(e,g) ~ | ~ g \in v_T^{-1}(T',e) ~\}, \end{align*}
and let $F(T)$ be a subgraph of $T$ with the set of vertices $V(F) = \bigcup_{(T',e) \in V(\G_T)} g_{T'}$, that is, $F(T)$ is a maximal subgraph of a fundamental domain of $v_T$.

\begin{prop}\label{prop-meager}
The set $C$ of compact leaves is a meager subset of $\fM_n$. Moreover, $C \cap \imath(X)$ is countable.
\end{prop}

\proof  Notice that $\cC = \imath^{-1}(C \cap \imath(X))$, where $\cC$ is defined as in the proof of Lemma \ref{prop-compactdense}. Then $C$ is meager if and only if $\cC$ is meager.

A leaf $L_T$ is compact if and only if $\cR(T)$ is a finite set. In this case $F(T)$ is compact. For each $r \in \mN $ define
  \begin{align*} A_r & = \{(T,e) \in X ~ | ~ F(T) \subset D_T(e,r) ~\}. \end{align*}
Since $D_{\cF_n}(e,r)$ contains at most a finite number of distinct subgraphs, there is a finite number of distinct $F(T) \subset D_T(e,r)$. Thus the set $A_r$ is finite and, therefore, nowhere dense in $X$. Since $\cC =  \bigcup_{r > 0} A_r $, the set $\cC$ is meager and countable.
\endproof

\begin{remark}
{\rm 
Alternatively, one can argue that given a compact leaf, there is always a one-ended leaf accumulating on the compact one, which shows that the set $C$ contains only leaves with non-trivial holonomy. Since the set of leaves without holonomy in every foliated space is residual by the result of Epstein, Millet, and Tischler \cite{EMT1977} (and also observed independently by Hector \cite{Hector2}) the set $C$ must be meager.
}
\end{remark}

\subsection{Fusion of graph matchbox manifolds}\label{section-directed}

In this section we prove Theorem \ref{thm-globalleaf}.(2), that is, that $(\cS_n,\preceq)$ is a partially ordered set, by means of the fusion technique, which associates to any given graph matchbox manifolds $\cM_1$ and $\cM_2$ a graph matchbox manifold $\cM_3$ such that $\cM_3 \supseteq \cM_1 \cup \cM_2$. 

Proof of Theorem \ref{thm-globalleaf}.(2) describes the fusion for a general case. An example of fusion is given in Section \ref{examples-prel}.

\proof \emph{(of Theorem \ref{thm-globalleaf}.(2)).} By Lemma \ref{defn-graphMM} it is enough to show that if $(T,e),(T',e) \in X$, then there exists $(S,e)$ such that $(T,e),(T',e) \in \overline{\cR(S)}$. The procedure of constructing $S$ consists of taking a $4$-ended graph $C_1$ and attaching to it copies of $D_T(e,r)$ and $D_{T'}(e,r)$, $r>0$, in a certain order, so that the graph $\G_S$ of the orbit $\cR(S)$ has $4$ ends, there is at least one end accumulating on $(T,e)$, and at least one end accumulating on $(T',e)$. Recall from Definition \ref{defn-holgraph} that the set of vertices $V(\Gamma_S) = \cR(S)$.

The \emph{leaf topology} on $\cR(S)$ is the topology induced from $\G_S$. Recall that a point $(T,e) \in \overline{\cR(S)}$ is in the limit set $\lim_{\bf e} \Gamma_S$ of an end ${\bf e} \in \cE(\G_S)$ if and only if there exists a sequence $\{(S_k,e)\} \in \Gamma_S$ converging to ${\bf e} \in \cE(\G_S)$ in the leaf topology, and to $(T,e)$ in the ambient topology. 

Choose $h,\wth \in F_n^0$, and let $C_0 \subset \cF_n$ be a subgraph with the set of vertices $V(C_0) = \bigcup_{n \in \mZ} \{h^n,\wth^n\}$. Then $C_{0}$ has $4$ ends. For $v \in \{h,\wth\}$ denote
  \begin{align*} {\bf e}_v^+ & = \lim_{n \to \infty} v^n, && {\bf e}_v^- = \lim_{n \to \infty} v^{- n}, \end{align*}
that is, $\{h^n\}_{n \in \mN}$ represents an end ${\bf e}_h^+$ and $\{h^{-n}\}_{n \in \mN}$ represents an end ${\bf e}_h^-$, and similarly for $v = \wth$. Let $c_v^+$ and $c_v^-$ respectively be infinite edge paths starting at $e$ and comprised from the shortest paths between the elements of the sets $\{v^n\}$, and $\{v^{-n}\}$ so that $C_{0} = c_h^+ \cup c_h^- \cup c_{\wth}^+  \cup c_{\wth}^-$. By Lemma \ref{lemma-projection} the graph $\G_{C_{0}}$ has four ends. We will obtain the graph $S$ by attaching finite decorations to $c_h^+$ and $c_{\wth}^+$ so that 
 \begin{align}\label{eq-fcond}\G_{T} \subset \lim_{{\bf e}_h^+} \G_{S} \cup \lim_{{\bf e}_{\wth}^+} \G_{S},\end{align}
and to  $c_h^-$ and $c_{\wth}^-$ so that 
  \begin{align}\label{eq-scond}\G_{T'} \subset \lim_{{\bf e}_h^-} \G_{S} \cup \lim_{{\bf e}_{\wth}^-} \G_{S}.\end{align}

If $T \ne T'$ there exists $r_0> 0 $ such that $D_T(e,r_0) \ne D_{T'}(e,r_0)$. If $T = T'$ then set $r_0 = 2$. Since $T$ and $T'$ are infinite, each of them has at least one end. Choose infinite edge paths $\delta$ in $T$, and $\delta'$ in $T'$ without intersections, such that $V(\delta)$ and $V(\delta')$ are sequences of points in $T$ and $T'$ converging to their ends.

For $i>0$ let 
  $$r_i = 2^{r_0+i} -1 ~ \textrm{and} ~ R_i = \sum_{k=0}^{i-1} r_k.$$ 
Let $w_i$ be an edge of $\delta$ such that $w_i \subset D_T(e,r_i)$ and $w_i \cap \partial D_T(e,r_i) \ne \emptyset$. Set 
  $$\beta = \left\{\begin{array}{ll} \phantom{-}1, & {\rm if} ~  t(w_i) \in \partial D_T(e,r_i), \\ -1, & {\rm if} ~ s(w_i) \in \partial D_T(e,r_i).\end{array} \right.$$
Let $v_i \in G_0$ be a letter labeling $w_i$. The intersection $\delta_i = D_T(e,r_i) \cap \delta$ is a finite path in $T$ of length $k_i$, and there is a finite word $E_i = h_1 \cdots h_{k_i-1}v_i^\beta$, $h_s \in G_0 \cup G_0^{-1}$,
such that $g \in \delta_i$ is a vertex in $\delta_i$ if and only if $g = h_1 \cdots h_{s}$ for $s < k_i$ or $g \in w_i \cap D_T(e,r_i)$. Set
   $$ \widetilde{E}_i = \left\{\begin{array}{ll} {\wth}^{R_i +r_i} E_i^{-1}, & {\rm if}~ v_i \ne \wth, \\ h^{R_i +r_i} E_i^{-1}, & {\rm if}~ v_i = \wth.\end{array} \right.$$
Let $C_i$ be a subgraph of $\cF_n$ containing $\widetilde{E}_i$ and such that there is an isomorphism $\alpha_i: C_i \to D_T(e,r_i)$ with $\alpha_i(\widetilde{E}_i) = e$. Then $\alpha_i(w_i) \cap C_0 \ne \emptyset$, and the union
  \begin{align*} S' = \bigcup_{i \geq 0} C_i \end{align*}
is a connected subgraph of $\cF_n$. By construction 
   \begin{align*} (T,e) \in \lim_{{\bf e}_h^+} \G_{S'} \cup \lim_{{\bf e}_{\wth}^+} \G_{S'}. \end{align*}
Implementing a similar algorithm for $(T',e)$, that is, attaching balls isomorphic to $D_{T'}(e,r_i)$, $r_i>r_0$, to the ends $e_{h}^-$ and $e_{\wth}^-$ of $S'$, one obtains a graph $S$ satisfying conditions \eqref{eq-fcond} and \eqref{eq-scond}.
\endproof

\subsection{Maximal transitive components} We prove Theorem \ref{thm-globalleaf}.(3), which says that $\fM_n$ contains recurrent leaves and so is a maximal element of $(\cS_n,\preceq)$. We need \cite[Theorem 3.5]{Blanc2001}, whose proof we give now for convenience of the reader.

\begin{lemma}\label{lemma-maximal}
Let $\cS' = \{\cM_i\}_{i \in \mN}$ be a totally ordered infinite subset of $(\cS_n,\preceq)$. Then 
  \begin{align*} \cM & = \overline{\bigcup_{i \in \mN} \cM_i}  \end{align*}
contains a recurrent leaf. Therefore, $\cS_n$ contains a maximal element.
\end{lemma}

\proof \cite[Theorem 3.5]{Blanc2001} Without loss of generality we can assume that every $\cM_{i-1}$ is a proper subset of $\cM_{i}$, and so no leaf $L \subset \cM_i$ is dense in $\cM$. Let $\{x_i\}$ be a sequence of points such that $x_i \in \cM_{i} \setminus \cM_{i-1}$. Since $\cM$ is closed it is compact and $\{x_i\}$ has a limit point $x$. By construction $x \notin \bigcup_{i \in \mN} \cM_i$, so $\cM \setminus \bigcup_{i \in \mN} \cM_i$ is non-empty.

Let $\{U_i\} \subset \cM$ be a system of open neighborhoods of $x$, such that $\bigcap_i U_i = x$. For every $U_i$ denote by $\widetilde{U}_i$ its saturation. The set $\widetilde{U}_i$ is open, and we claim that $\widetilde{U}_i$ is dense in $\cM$.

Indeed, let $V \subset \cM$ be open and let $\widetilde{V}$ be its saturation. Then there exists $\cM_{k_i}, \cM_{k_V} \in \cS'$ such that $\widetilde{U}_i \cap \cM_{k_i} \ne \emptyset$, and $\widetilde{V} \cap \cM_{k_V} \ne \emptyset$. Since $\cS'$ is totally ordered, we have either $\cM_{k_i} \preceq \cM_{k_V}$ or $\cM_{k_V} \preceq \cM_{k_i}$. For definitiveness assume the former. Let $L \subset \cM_{k_V}$ be a transitive leaf, then $L \subset \widetilde{U}_i \cap \widetilde{V}$, which implies that $\widetilde{U}_i$ is dense in $\cM$.

The intersection $\bigcap_i \widetilde{U}_i$ of a countable family of open dense subsets is exactly the leaf $L_x \owns x$. Since $\cM$ is compact and Hausdorff, it is a Baire space, and $L_x$ is dense in $\cM$. Since $L_x$ is in the boundary of $\overline{\bigcup_{i \in \mN} \cM_i}$, $L_x$ is recurrent and $\cM \in \cS_n$.

It follows that every totally ordered subset of $(\cS_n,\preceq)$ has an upper bound which is a recurrent graph matchbox manifold if $\cS'$ is an infinite chain, and a transitive graph matchbox manifold in the case when $\cS'$ is finite. Applying Zorn's lemma we conclude that $\cS_n$ contains a maximal element.
\endproof

Recall \cite{Kur} that a set $A \subset M$ is \emph{residual} if its complement $M \setminus A$ is meager.

\proof \emph{(of Theorem \ref{thm-globalleaf}.(3))} Let $\cM_{max}$ be a maximal element in $(\cS_n, \preceq)$, and $\cM$ be any matchbox manifold. By Theorem \ref{thm-globalleaf}.(2) there exists a matchbox manifold $\cM'$ obtained by fusion, such that $\cM_{max} \cup \cM$ is a subset of $\cM'$. Since $\cM_{max}$ is maximal, $\cM_{max} = \cM'$ and, therefore, $\cM \subset \cM_{max}$. It follows that there is a unique maximal element which contains every graph matchbox manifold $\cM$. Then necessarily $\cM_{max} = \fM_n \in (\cS_n,\preceq)$. Since $X_n$ is a Cantor set, by Corollary \ref{recurrencelemma} $\fM_n$ is recurrent.
\endproof

\section{Theory of levels for graph matchbox manifolds} \label{exmpl} 

In this section, inspired by ideas of Cantwell and Conlon \cite{CC1, CC3}, and also of Hector \cite{Hector}, Nishimori \cite{Nish1977,Nish1979} and Tsuchiya \cite{Ts1980} for codimension $1$ transversally $C^2$-differentiable foliations, we introduce a quantifier of dynamical complexity of a set in $(\cS_n,\preceq)$ which, following \cite{CC1}, we call \emph{level}. We repeat now Definition \ref{defn-level}.

\begin{defn}\label{defn-levels}
Let $\cM \subset \fM_n$ be a graph matchbox manifold. 
\begin{enumerate}
\item $\cM$ is said to be \emph{at level $0$} if either $\cM$ is a compact leaf, or $\cM$ is a minimal foliated space. In that case all leaves of $\cM$ are also at level $0$. 
\item $\cM$ is \emph{at level $k$} if the closure of the union of leaves which are not dense in $\cM$, is a proper closed subset of $\cM$, every such leaf is at level at most $k-1$, and there is at least one leaf at level $k-1$. If $\cM$ is at level $k$, a leaf $L \subset \cM$ is at level $k$ if and only if $L$ is dense in $\cM$.
\item A leaf $L$ is at \emph{infinite level} if it is not at finite level.
\end{enumerate}
\end{defn}

The goal of this section is to prove Theorems \ref{thm-infiniteleaves}, \ref{thm-examples1totprop} and \ref{thm-examples1recurrent}. 

\subsection{Leaves at infinite levels}

We investigate leaves at infinite levels and prove Theorem \ref{thm-infiniteleaves}. We first notice that Definition \ref{defn-levels} allows for two types of leaves at infinite levels. 

Let $\cM = \overline{L}$. Then $\overline{L} \subset \fM_n$ is not at a finite level if one of the following holds:

\begin{enumerate}
\item {\bf Type 1}. The union of leaves which are not dense in $\cM$ is dense in $\cM$, and $\cM$ is a proper subset of $\fM_n$. 
\item {\bf Type 2}. The union of leaves which are not dense in $\cM$ is a proper subset of $\cM$ and contains a leaf at infinite level.
\end{enumerate}

Theorem \ref{thm-infiniteleaves} states that $\fM_n$ contains leaves of both types. A proof of that is given in the following two propositions.

\begin{prop}
In the space of graph matchbox manifolds $\fM_2$ there exists a leaf $L$ and a graph matchbox manifold $\cM = \overline{L}$ at infinite level of Type $1$.
\end{prop}

\proof We have to prove that there exists an infinite increasing chain
  \begin{align*} \cM_0 \subset \cM_1 \subset \cM_2 \subset \cdots\end{align*}
of distinct graph matchbox manifolds such that $\cM = \overline{\bigcup_i \cM_i}$ is a proper subset of $\fM_n$. Then $\cM$ is a graph matchbox manifold at infinite level of Type 1.

Let $F_2^0 = \{a,b\}$, and let $T_0$ be a subgraph of $\cF_2$ with the vertex set $V(T_0) = \cup_{n \geq 0}\{ a^n,a^{-n}\}$. then $L_0 = L_{T_0}$ is a compact leaf. Let $T_1$ be a subgraph of $\cF_2$ with the vertex set $V(T_1) = \cup_{n \geq 0} a^n$. Then $L_1 = L_{T_1}$ is a proper leaf, and $\cM_1 = \overline{L_1}$ contains $L_1$ and $L_0$. For $n>2$ construct a graph $T_n$ by applying the fusion construction of Section \ref{section-directed} for $T = T'=T_i$. Define the \emph{depth} $k(T_n)$ to be the maximal number of changes from an edge labeled by $a$ to an edge labeled by $b$ or vice versa in a finite path in $T_n$ starting at the origin. It is easy to see that $k(T_n) \geq k(T_{n-1}) +1$. This implies that for $R>0$ large enough to contain a path of depth $k(T_n)$, there is no isomorphism from a closed ball in $T_i$, $i< n$, onto a closed ball of radius $R$ around the origin in $e$, and so $T_n \notin \bigcup_{0 \leq i <n} \cR(T_i)$. Therefore, all $T_n$ and $\cM_n = \overline{L_{T_n}}$ are distinct and the complement of $\cup_{n \geq 0} \cM_n$ in the closure
$$\cM = \overline{\bigcup_{n \geq 0} \cM_n}$$
is non-empty. Then by Lemma \ref{lemma-maximal} $\cM$ is a recurrent graph matchbox manifold at infinite level of Type $1$. We notice that $\cM$ is a proper subset of $\fM_2$ since it is easy to find a finite subgraph $S \subset B_{\cF_2}(e,3)$ which is isomorphic to no subgraph of $T_n$ for every $n \geq 0$. 
\endproof

\begin{prop}\label{infinite-hierarchy}
Let $\cM$ be a graph matchbox manifold at infinite level and suppose $\cM$ is a proper subset of $\fM_n$. Then there exists a graph matchbox manifold $\widetilde{\cM}$ such that 
  $$\cM \subset \widetilde{\cM} \subset \fM_n$$ 
are proper inclusions. The space $\widetilde{\cM}$ is at infinite level of Type 2.
\end{prop}

\proof Let $L \subset \cM$ be a dense leaf, and $(T_L,e) \in \imath^{-1}(L \cap \imath(X))$. By assumption the complement $\fM_n \setminus \cM$ is open, then we can find a compact leaf $C$ such that $C \cap \cM = \emptyset$. Let $T_C \in \imath^{-1}(C \cap \imath(X))$, and let $T_\Phi$ be a graph obtained by fusion on $T_L$ and $T_C$ (see Section \ref{section-directed}), and $\Phi$ be the corresponding leaf. Then $\overline{\Phi} \supset \overline{L} \cup C$. The complement of $\overline{\Phi}$ in $\fM_n$ is open, and to prove the proposition we have to show that $\Phi$ can be chosen in such a way that the complement of $\overline{\Phi}$ is non-empty.

The leaf $\Phi$ can either be recurrent or proper. If $\Phi$ is proper, then every point of $\cR(\Phi)$ is isolated in $\overline{\cR(\Phi)}$, and since $X_n$ is a Cantor set, this implies that the complement $\fM_n \setminus \overline{\Phi}$ is non-empty. We show that $\Phi$ is proper.

Since $C \cap \cM = \emptyset$, and $\overline{\Phi} \supset C$, $\overline{L}$ is a proper subset of $\overline{\Phi}$. Suppose $\Phi$ accumulates on itself, then by Lemma \ref{lemma-recurrencecritgraph} there exists a sequence $\{g_\ell\} \in V(T_\Phi)$, $g_\ell \ne e$, such that for every $r>0$ there exists $\ell_r>0$, and for all $\ell \geq \ell_r$ there is an isomorphism 
  $$\alpha^r_\ell: D_{T_\Phi}(e,r) \to D_{T_\Phi}(g_\ell,r).$$
In the fusion construction of Section \ref{section-directed} let $h, \wth \in F_n^0$ and let $A$ and $B$ be subgraphs of $\cF_n$ with $V(A) = \cup_{n \in \mZ} h^n$ and $V(B) = \cup_{n \in \mZ} \wth^n$. Recall that $T_{\Phi}$ was obtained by attaching copies of closed balls $D_{T_L}(e,n)$ and $D_{T_C}(e,n)$ to the union of $A$ and $B$ in $\cF_n$. Without loss of generality we can assume that $D_{T_\Phi}(g_{\ell},r)$ is not isomorphic to a subset of $D_{T_L}(e,n)$ for any $\ell,r,n>0$. Indeed, suppose there exists a subsequence $\{g_{\ell_k}\}$ such that for $r>r_k$ there is an isomorphism onto its image
   \begin{align}\label{eq:choiceofballs}D_{T_\Phi}(g_{\ell_k},r) \to D_{T_L}(e,n_k),\end{align}
$n_k\geq r$. Then $\Phi \subset \overline{L}$ which contradicts the construction of $\Phi$. By a similar argument we can assume that $D_{T_\Phi}(g_{\ell},r)$ is not isomorphic to a subset of $D_{T_C}(e,n)$ for any $\ell,r,n>0$.

Suppose $g_\ell$ is a vertex in the copy of $D_{T_L}(e,n_\ell)$ attached to $A \cup B$, and let $a_\ell$ be the length of the shortest path between $g_\ell$ and the point of the attachment. Then either $\{a_\ell\}$ contains a monotonically increasing subsequence, or $\{a_\ell\}$ converges to an integer $a$. If the first situation occurs, then it is possible to choose a subsequence of $\{g_\ell\}$ satisfying the property \eqref{eq:choiceofballs}, which contradicts the choice of $\Phi$. The second situation cannot occur either by the following argument. Without loss of generality we can assume that $g_\ell$ is at distance $a$ from $A$, and $v_\ell \in A$ is such that $d_{T_\Phi}(g_\ell,v_\ell) = a$. Then for $r>a$ the image of $v_\ell$ is either in $B$, or inside one of the attached copies of $D_{T_L}(e,r)$. The first situation is not possible, since by construction the image of $v_\ell$ has two adjacent $a$-edges, while every point in $B$ except the origin has at most one adjacent $a$-edge. In the second case the isomorphism between $D_{T_\Phi}(g_\ell,r)$ and $D_{T_\Phi}(e,r)$, $r>a$, implies that the copy of $D_{T_L}(e,r)$ in $\Phi$ which contains the image of $v_\ell$, also contains an infinite subgraph, which is not possible. By a similar argument $\{g_\ell\}$ cannot contain a subsequence of vertices lying in the copies of $D_{T_C}(e,n)$. 

Therefore, $\{g_\ell\} \in V(A \cup B)$. Passing to a subsequence without loss of generality we can assume that $g_\ell \in V(A)$. But such a sequence cannot exist since every point of $A$ except the origin is a vertex with at most three adjacent edges, while the origin is a vertex with four adjacent edges. It follows that $\Phi$ is a proper leaf, and the complement $\fM_n \setminus \overline{\Phi}$ is non-empty. Then $\widetilde{M} = \overline{\Phi}$ is at infinite level of Type $2$.
\endproof

\subsection{Proof of Theorem \ref{thm-examples1totprop}}\label{sec-totprop} Recall \cite{CandelConlon2000} that a graph matchbox manifold $\cM$ is \emph{totally proper} if it contains only proper leaves \cite{CC1}. Theorem \ref{thm-examples1totprop} states that there exists a totally proper graph matchbox manifold $\cM$ at level $2$ which contains a countably infinite number of distinct compact leaves. First we need a few auxiliary lemmas. Let $F_n^0 = \{a,b\}$.

\begin{figure}
\begin{minipage}{0.5\linewidth}
\centering
\includegraphics [width=3.5cm] {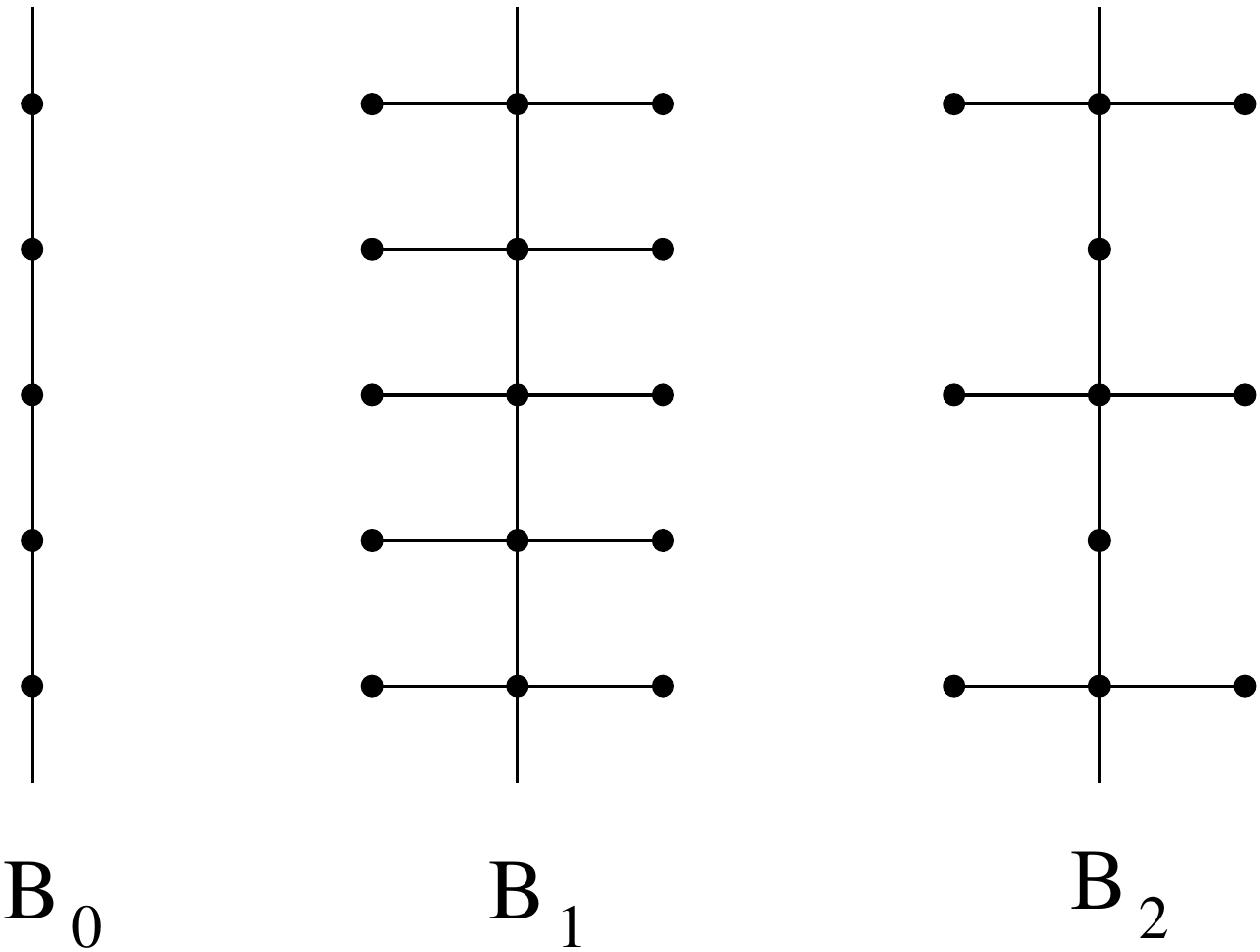}
\caption{Graphs $B_m$, $m=0,1,2$.}
 \label{fig:graphsBn}
\end{minipage}%
\begin{minipage}{0.5\linewidth}
\centering
\includegraphics [width=2.5cm] {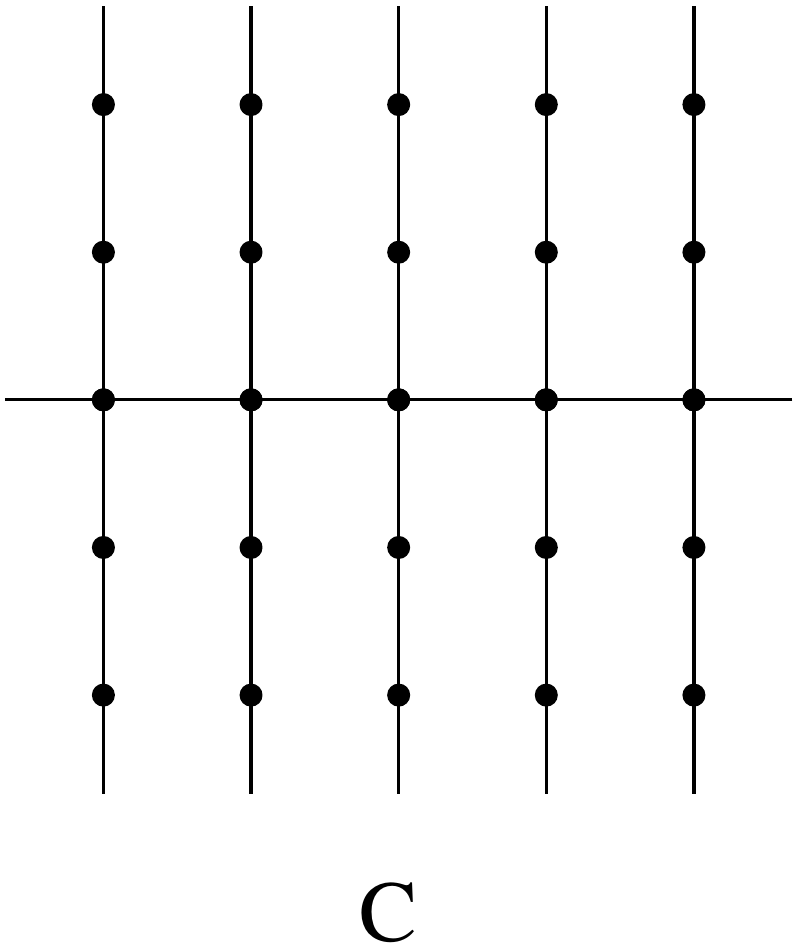}
\caption{Graph $C$.}
 \label{fig:graphC}
\end{minipage}
\end{figure}

\begin{lemma}\label{compactleaves}
Given $m \geq 0$, there exists $(B_m,e) \in X$ such that the leaf $L_m = L_{B_m} \subset \fM_n$ is compact, and all leaves $L_m$ are distinct.
\end{lemma}

\proof 
Let $B_0$ be a subgraph of $\cF_2$ with the set of vertices $V(B_0) =\bigcup_{k \geq 0} \{b^k,b^{-k}\}$ (see Fig.\ref{fig:graphsBn}). Given $m \in \mN$ let $B_m$ be a subgraph of $\cF_2$ with the set of vertices
  \begin{align*} V(B_m) & = V(B_0) \bigcup \left( \bigcup_{k \in \mZ} \{ b^{mk} a, b^{mk} a^{-1}\} \right).  \end{align*}
Then $B_m$ is invariant under the action of $b^{mk}$, and the corresponding leaf $L_m$ is compact and homeomorphic but not isometric to the standard $2$-torus. \endproof

\begin{lemma}\label{leaflevel1}
There exists a graph $(C,e) \in X$ such that the leaf $L_C \subset \fM_n$ contains a single leaf $L_0$ in its limit set. It follows that $\cM_C$ is at level $1$.
\end{lemma}

\proof  Let $A_0$ be a subgraph of $\cF_2$ with the set of vertices $V(A_0) =\cup_{k \geq 0} \{a^k,a^{-k}\}$. For every $k > 0$ let $\widetilde{B}_{\delta k}$, $\delta = \pm 1$, be a subgraph of $\cF_2$ containing $a^{\delta k}$ such that there is an isometry $\alpha_{\delta k}: B_0 \to \widetilde{B}_{\delta k}$ with $\alpha_k(e) = a^{\delta k}$. Set (see Fig.\ref{fig:graphC})
  \begin{align*} C & = A_0 \cup B_0 \cup \left( \bigcup_{k > 0} \{\widetilde{B}_k, \widetilde{B}_{-k}\} \right). \end{align*}
The pseudogroup $\fG_C$ of the corresponding matchbox manifold $\cM_{C}$ is generated by $\gamma_a$ and $\gamma_b$, and the graph $(C,e)$ is invariant under $\gamma_a$. We have
  \begin{align*} {\rm dom}(\gamma_a) \cap \cR(C) & = (C,e), && \cR(C) \subset {\rm dom} (\gamma_b).  \end{align*}
Clearly $(B_0,e) \subset \overline{\cR(C)}$, so $L_0 \subset \lim L_C$. We now have to show that there are no other leaves in $\overline{L_C}$.

Suppose $(T,e) \in \overline{\cR(C)}$. Then there exists a sequence $(C,g_s)$, $s \in \mN$, such that for every $r> 0 $ there exists $s_r$ such that for all $s>s_r$ there is an isometry
   $$\alpha_s^r:  D_C(g_s,r) \to D_T(e,r).$$ 
Since $\gamma_a(C,e) = (C,e)$ we can assume that $g_s = b^{k_s}$, $k_s \in \mZ$. Suppose there exists $m \in \mN$ such that for $s>m$ and all $r>0$ we have $D_C(b^{k_s},r-1)  \cap A_0 = \emptyset$. Then $(T,e) = (B_0,e)$. If not, let $m>0$ be such that $ D_C(b^{k_m},r-1) \cap A_0 \ne \emptyset$. Then then for all $s> m$ we have $d_C(b^{k_s},e) = d_C(b^{k_m},e)$ and $b^{k_s} = b^{k_m}$. Then $(T,e) = (C,b^{k_m})$. Thus $\cM_C$ contains only $2$ leaves, one of which is at level $0$, and it follows that $\cM_{C}$ is at level $1$.
\endproof

\begin{remark}\label{remark-anotherexample}
{\rm
Of course, there are many examples which satisfy requirements of Lemma \ref{leaflevel1}. Let $B_0'$ be a subgraph of $\cF_2$ with the set of vertices $V(B_0') = V(B_0) \cup \{a,a^{-1}\}$. Then $L_0 \subset \cM_{B_0'}$, and an argument similar to the one in Lemma \ref{leaflevel1} shows that $\cM_{B_0'}$ is at level $1$.
}
\end{remark}

The following proposition completes the proof of Theorem \ref{thm-examples1totprop}.

\begin{prop}\label{leafinfiniteends}
There exists $(T,e) \in X$ such that the graph matchbox manifold $\cM_{T}$ has the following properties.
\begin{enumerate}  
\item $\bigcup_m L_m  \subset \cM_{T}$, where $L_m$ are as in Lemma \ref{compactleaves}.
\item For every $m \in \mN$ there is a distinct end ${\bf e}_m$ such that $L_m \subset \lim_{{\bf e}_m} L_T$.
\item $\cM_{T}$ is at level $2$.
\end{enumerate}
\end{prop}

\proof Let $A_0$ be a subgraph of $\cF_2$ with the set of vertices $V(A_0) =\cup_{k \geq 0} \{a^k,a^{-k}\}$. For every $k > 0$ let $\widetilde{B}_{\delta k}$, $\delta = \pm 1$, be a subgraph of $\cF_2$ containing $a^{\delta k}$ such that there is an isometry $\alpha_{\delta k}: B_k \to \widetilde{B}_{\delta k}$ with $\alpha_{\delta k}(e) = a^{\delta k}$, where $B_k$ was constructed in Lemma \ref{compactleaves}. Set 
  \begin{align*} T & = A_0 \cup B_0 \cup \left( \bigcup_{k > 0} \{\widetilde{B}_k, \widetilde{B}_{-k}\} \right). \end{align*}
Then $(1)$ and $(2)$ are clearly satisfied. We have to show that $\cM_{T}$ is at level $2$.

We first notice that $\cM_T$ is at level at least $2$, since $\cM_{T} \supset \cM_{C}$ where $\cM_{C}$ was obtained in the proof of Lemma \ref{leaflevel1}. To see that notice that, given $r>0$, for every $s \geq s_r = r$ there is an isometry $\alpha_s^r: D_T(a^{2s},r) \to D_C(e,r)$, and so the sequence $\{(T,e) \cdot a^{2r}\}$ converges to $(C,e)$. 

We now show that non-dense leaves in $\cM_T$ are at level at most $1$, and so $\cM_T$ is at level $2$. Suppose $(T',e) \in \overline{\cR(T)}$. Then there is a sequence $(T,g_s)$ converging to $(T',e)$, i.e. for every $r\geq 2$ there exists $s_r>0$ such that for every $s \geq s_r$ there is an isometry $\alpha_s^r: D_T(g_s,r) \to D_{T'}(e,r)$. The following argument is just a more complicated version of the one in Lemma \ref{leaflevel1}. Let
 \begin{align*}\ell & = \min_s \{~{\rm dist} (g_s,A_0) ~| ~ s \in \mN ~\}, \end{align*}
and let $s_1$ be such that ${\rm dist}(g_{s_1},A_0) = \ell$. Consider the following cases.

\emph{Case 1.} Suppose $\ell = 0$.  Then $g_{s_1} = a^{m_1}$, $m_1 \in \mZ$. Since $r \geq 2$, $D_T(g_{r_1},r)$ contains an edge path of $a$-edges of length at least $4$, and it follows that for all $r\geq r_1$ and all $s \geq s_1$ we have 
  $$g_{s} = a^{m_s}, ~m_s \in \mZ.$$ 
Then two situations are possible: either there exists $r_2 \geq r_1$ and $s_2 \geq s_1$ such that $D_T(g_{s_2},r_2)$ contains another path of $a$-edges which necessarily has length $2$, or such an $r_2$ does not exist. In the first case it follows that for all $r \geq r_2$ and all $s \geq s_2$ we have $g_s = a^{\pm m_2}$, $m_2 \in \mZ$. Since $\{(T,e) \cdot g_s\}$ converges, there is an $s_3 \geq s_2$ such that for all $s \geq s_3$ either $g_s = g_{s_3} = a^{m_2}$ or $g_s = g_{s_3} = a^{-m_2}$. If the former is true, then $(T',e) = (T,e) \cdot a^{m_2}$, otherwise $(T',e) = (T,e) \cdot a^{-m_2}$. In the second case for every $r \geq r_1$ and $s \geq s_1$ one constructs an isomorphism $D_T(g_s, r) \to D_C(e,r)$, and it follows that $(T',e) = (C,e)$.

\emph{Case 2.} Suppose $\ell> 0$, and set $\ell_s = {\rm dist}(g_s,A_0)$, $s \in \mN$. Let $m_s$ be such that $\ell_s = d_T(a^{m_s},g_s)$.

\emph{Case 2.1.} Suppose there exists $r_1\geq 2$ and $s_1 >0$ such that for all $r \geq r_1$ and all $s \geq s_1$ we have $\ell_s \geq r$. Then either
 $$g_s = a^{m_s}b^k, ~ |k|>1,$$ 
or 
 $$g_s = a^{m_s}b^k a^\beta, ~ |k|>1, ~ \beta = \pm 1.$$ 
Then three situations are possible. 

\emph{Case 2.1.1.} First, there may exist $r_2 \geq r_1$ and $s_2 \geq s_1$  such that the ball $D_T(g_{s_2},r_2)$ contains at least two distinct paths of $a$-edges $c_1$ and $c_2$, which are necessarily of length $2$. Let
  $$\tilde{\ell}  = d_H(c_1,c_2),$$
where $d_H$ is the distance between these sets. Then for every $r>r_2$ the ball $D_T(g_{s},r)$ must contain two distinct paths of $a$-edges such that the distance between them is exactly $\tilde{\ell}$. Then either $g_{s_2} \in V(\widetilde{B}_{\tilde{\ell}})$ or $g_{s_2} \in V(\widetilde{B}_{-\tilde{\ell}})$, and the same is true for $g_s$ for $s \geq s_2$. It follows that $(T',e) \in \cR(B_{\tilde{\ell}})$, as we show now.

Indeed, if $g_s = a^{m_s}b^k$, let $\tilde{g} \in V(T)$ be such that
  \begin{align*} d_T(g_s,g_s\tilde{g})  & = {\rm dist}(g_s,c_1 \cup c_2), \end{align*}
and if $g_s = a^{m_s}b^k a^\beta$, let $\tilde{g} = a^{- \beta}$. Then for every $r \geq r_2$ and every $s \geq s_2$ there is an isometry
 $$\alpha^r_s : D_T(g_s,r) \to D_{B_{\tilde{\ell}}}(\tilde{g}^{-1},r),$$ 
which implies $(T',e) = (B_{\tilde{\ell}}, e) \cdot \tilde{g}^{-1}$. 

\emph{Case 2.1.2.} Second, there may exist $r_2 \geq r_1$ and $s_2 \geq s_1$ such that for all $r\geq r_2$ and all $s \geq s_2$ the ball $D_T(g_{s},r)$ contains exactly one path of $a$-edges $c_1$, which is necessarily of length $2$. Then $r' > r$ implies $s_{r'} > s_r$, that is, there is a subsequence $\{g_{s_v}\} \subset \{g_s\}$ such that distinct $g_{s_v}$ lie in distinct $\widetilde{B}_{\tilde{l}_v}$. Then by a similar argument as in \emph{Case 2.1.1} $(T',e) \in \cR(B_0')$, where $B_0'$ is as in Remark \ref{remark-anotherexample}. 

\emph{Case 2.1.3.} Third, if for all $r \geq r_1$ and all $s \geq s_1$ the closed ball $D_T(g_s,r)$ does not contain any $a$-edges, then there is an isometry 
  $$D_T(g_s,r) \to D_{B_0}(e,r),$$ 
and it follows that $(T',e) = (B_0,e)$.

\emph{Case 2.2.} Suppose there exists $r' \geq 2$ and $s' > 0$ such that for all $r \geq r'$ and all $s > s'$ we have $\ell_s < r$. Then $D_T(g_s,r) \cap A_0 \ne \emptyset$ and two situations are possible. 

\emph{Case 2.2.1.} There exists $r_1 \geq 2$ and $s_1 >0$ such that for all $r \geq r_1$ and all $s \geq s_1$ we have $\ell_s < r-1$. Then $D_T(g_s,r)$ contains a path of $a$-edges of length at least $4$, and for $s>s_1$ we have $\ell_s = \ell_{s_1}$. If there exists $r_2 \geq r_1$ and $s_2 \geq s_1$ such that $D_T(g_{s_2},r)$ contains another path of $a$-edges, which is of length necessarily $2$, then for all $s\geq s_2$ either $g_s = g_{s_2}$ or $g_s = g_{s_2}^{-1}$. By increasing $r$ and possibly $s_2$ we exclude one of these options, and it follows that either $(T',e) = (T,e) \cdot g_{s_2}$ or $(T',e) = (T,e) \cdot g_{s_2}^{-1}$. If $r_2$ with this property does not exist, then $g_s = a^{m_s}b^k$ for some $k \in \mZ$, and $(T',e) = (C,e) \cdot b^k$.

\emph{Case 2.2.2.} If an $r_1$ as specified in \emph{Case 2.2.1} does not exist, then for all $r \geq r'$ and all $s>s'$ we have $\ell_s = r-1$, and $D_T(g_s,r)$ contains at least one path of $a$-edges of length $2$. If there exists $r_2 \geq r'$ and $s_2 > s'$ such that $D_T(g_{s_2},r_2)$ contains at least $2$ distinct $a$-edge paths, then we are in the situation of the \emph{Case 2.1.1}, and $(T',e) \in \cR(B_k)$ for some $k \geq 1$. Otherwise $D_T(g_{s},r)$ contains exactly one $a$-edge path of length $2$, and we are in the situation of the \emph{Case 2.1.2}. Then $(T',e) \in \cR(B_0')$.

\emph{Case 2.3.} If $r'>2$ and $s'>0$ such as in Case 2.1. and Case 2.2 do not exist, then there is a subsequence $\{g_{k_s}\}$ of $\{g_s\}$ such that there exists $r' \geq 2$ and $s' > 0$ such that for all $r \geq r'$ and all $s > s'$ we have $\ell_{k_s} < r$ (resp. $\ell_{k_s} \geq r$). Then use the argument of \emph{Case 2.2.} (resp. \emph{Case 2.1.}). 

This exhausts the distinct possibilities for choosing a sequence $\{g_s\}$. We have shown that $\lim L_T$ contains leaves at level at most $1$. It follows that $\cM_T$ is at level $2$.
\endproof

\subsection{Proof of Theorem \ref{thm-examples1recurrent}}\label{ssec-recurrent} We prove Theorem \ref{thm-examples1recurrent}, which states that there exists a graph matchbox manifold at level $2$ exhibiting interesting dynamics, namely, one can find a clopen subset of a transversal such that the restricted pseudogroup is equicontinuous. 

We construct such a matchbox manifold for $n=2$, i.e. in $\fM_2$. Let $F_n^0 = \{a,b\}$ be a set of generators. We first give a construction of a tree $(K,e) \in X$, and then prove that $\cM_K$ satisfies the requirements of the theorem in a series of propositions and lemmas.

Let $C_0= L_0$ be a subgraph of $\cF_2$ with the set of vertices $V(C_0) = V(L_0) = \{e,a,a^{-1},b,b^{-1}\}$. Given $C_i$ and $L_i$ define $A_i^{\delta}$, $\delta = \pm 1$, to be a subgraph of $\cF_2$ containing the vertex $a^{\delta 2^i}$ such that there is an isomorphism
  \begin{align*} \alpha_i^\delta & : A_i^\delta \to C_i ~ ~{\rm with}~ ~ \alpha_i^\delta(a^{\delta 2^i}) = e, \end{align*}
and $B_i^{\delta}$ be a subgraph of $\cF_2$ containing the vertex $b^{\delta 2^i}$ such that there is an isomorphism
  \begin{align*} \beta_i^\delta & : B_i^\delta \to L_i ~ ~{\rm with}~ ~ \beta_i^\delta(b^{\delta 2^i}) = e. \end{align*}

\begin{figure}
\centering
\includegraphics [width=5.8cm] {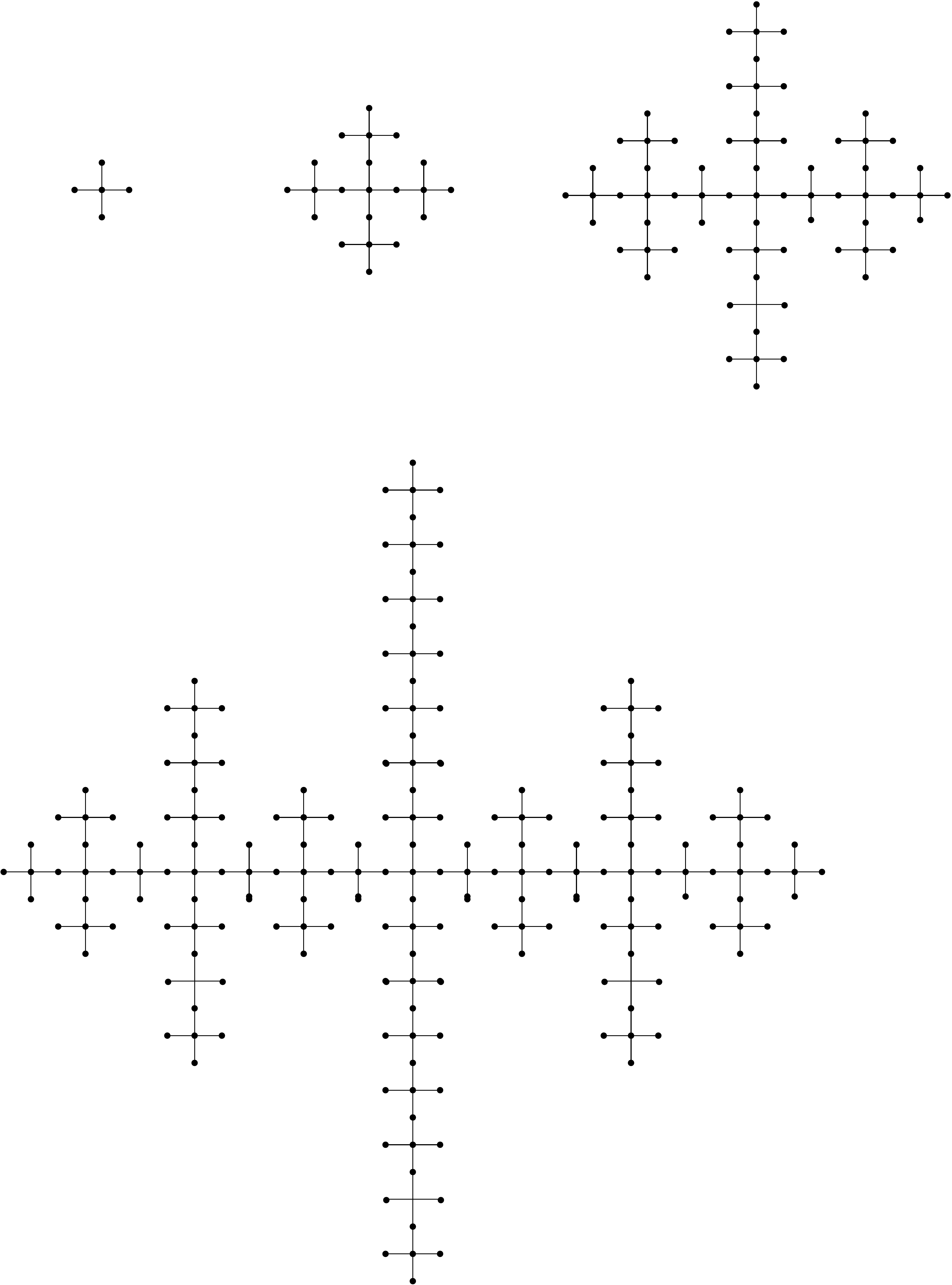}
\caption{Construction of $K$.}
 \label{fig:modif-kenyon}
\end{figure}

Then set $C_{i+1} = C_i \cup A_i^{1}\cup A_i^{-1} \cup B_i^1 \cup B_i^{-1}$, and $L_{i+1} = L_i \cup B_i^{1} \cup B_i^{-1}$. Define $K = \bigcup_{i \in \mN} C_i$ (see Figure \ref{fig:modif-kenyon}). We show that $\cM_{K}$ contains an uncountable number of leaves, and each leaf is at finite level.

\begin{lemma}\label{prop-recurrent}
The graph matchbox manifold $\cM_{{K}}$ is recurrent, and, therefore, it has an uncountable number of leaves.
\end{lemma}

\proof Let $g_i = a^{2^i}$, $i\in \mN$, so all $g_i$ are distinct. Then for $k \geq i$ there is an isometry 
  $$\alpha^i_k: D_{K}(g_k,2^i-1) \to D_{K}(e,2^i-1),$$ 
and by Lemma \ref{lemma-recurrencecritgraph} $\cM_{K}$ is recurrent. Therefore, $\cM_{{K}}$ contains an uncountable number of leaves.
\endproof

Let $A \subset K$ be a subgraph with the set of vertices $V(A) = \bigcup_{i \in \mZ} \{a^i\}$. We call a connected subgraph $F \subset K$ a \emph{vertical decoration} of $K$, if the set of vertices $V(F)$ satisfies the following two conditions. 
\begin{enumerate}
\item There exists a vertex $a^{n}$, $n \in \mZ$, called the \emph{point of attachment} of $F$ to $A$, such that $a^{n} \in V(F)$ and for every $m \ne n$ we have $a^{m} \notin V(F)$.
\item If $a^nb^k \in V(F)$ then $a^nb^{-k} \notin V(F)$.
\end{enumerate}

A vertical decoration $F$ is a one-ended connected subtree of $K$ with exactly one vertex $a^n$ lying in $A$. We say that $F$ is \emph{finite} of length $\ell$ if $\ell$ is the length of a longest path without self-intersections contained in $F$. If $F$ is not finite, we say that $F$ is \emph{infinite}. We notice that $K$ has $2$ distinct infinite decorations attached to $A$ at the vertex $e$, and all other vertical decorations are finite.

To prove Theorem \ref{thm-examples1recurrent} we need the following two lemmas.

\begin{lemma}\label{lemma:infinitedecorations}
If $(T,e) \in \overline{\cR(K)}$ has more than $2$ ends, then $(T,e) \in \cR(K)$.
\end{lemma}

Lemma \ref{lemma:infinitedecorations} states that $L_K$ is the only $4$-ended leaf in $\cM_K$, and every other leaf is either $2$-ended, $1$-ended or compact.

\proof Let $(T,e) \in \overline{\cR(K)}$, and notice that $(T,e)$ cannot have more than two infinite vertical decorations. Indeed, suppose there are more than two such decorations, then there are at least two distinct vertices $v_1$ and $v_2$ where these decorations are attached to the subgraph $A$. Since $(T,e) \in \overline{\cR(K)}$, there exists a sequence $\{g_s\} \in V(K)$ such that for every 
  $$r > \max \{d_T(e,v_1), d_T(e,v_2)\}$$  
there is $s_r>0$ such that for all $s>s_r$ there is an isometry $B_{K}(g_s,r) \to D_{T}(e,r) $, which implies that for any $m>0$ we should be able to find a pair $F_1^{(m)}, F_2^{(m)}$ of vertical decorations in $K$ of length at least $m$ such that the distance between their points of attachment to $A$ is precisely $d_T(v_1,v_2)$. This is not possible, since the distance between decorations of length at least $m$ in $K$ strictly increases with $m$. Therefore, $(T,e)$ has at most two infinite vertical decorations.
Now suppose $(T,e)$ has two infinite vertical decorations, and let $v_1 \in T$ be the vertex at which the decorations are attached. Then for any $r>0$ one can construct an isometry $D_K(e,r) \to D_T(v_1,r)$, which implies $(T,e) = (K,e) \cdot v^{-1}_1$. \endproof

To keep track of leaves in $\cM_K$, we construct a section of $\overline{\cR(K)} \setminus \cR(K)$ using the coding procedure, employed in \cite[Section 3.2]{ALM2008} to study the example of Kenyon and Ghys.

\begin{lemma}\label{leavesfinitelevel}
Let $\cQ_4 = \left\{ a,a^{-1},b,b^{-1} \right\}^\mN$ be the set of one-sided infinite sequences. Then there exists a subset $\cQ \subset \cQ_4$ and a map $P:\cQ \to X$ such that for every $(T,e) \in \overline{\cR(K)} \setminus \cR(K)$ we have $P(\cQ) \cap \cR(T) \ne \emptyset$.
\end{lemma}

\proof
Consider a subset 
  $$\cQ = \bigcup_{n \in \mN}\{(\alpha_0 \ldots \alpha_n (a a^{-1})) ~ | ~  \alpha_i = b^{\delta},~ \delta \in \{-1,1\}, ~ 0 \leq i \leq n \} \cup \{(bb^{-1})\} \cup\{(b)\}\cup \{(b^{-1})\}  \subset \cQ_4,$$ 
and obtain the map $P: \cQ \to X$ as follows. Let $x_0 = e$, and $x_1 = x_0 \alpha_0 \in V(\cF_2)$. Set $E_1 = L_1$ to be the subgraph of $\cF_2$ with the set of vertices $V(E_1) = V(L_1) = \{x_1, x_1a,x_1a^{-1},x_1b,x_1b^{-1}\}$. For $i \geq 2$ obtain the graphs $(E_i,x_i)$ are by the following inductive procedure. 

Let $x_{i-1}$, $L_{i-1}$ and $E_{i-1}$ be given. Let $x_{i} = x_{i-1}\alpha_{i-1}^{2^{i-1}}$, then there is an edge $w \subset \cF_2$ such that $E_{i-1} \cap w$ is a vertex, and such that $x_i$  is another vertex of $w$. Set $L'$ and $E'$ to denote the subgraphs of $\cF_2$ containing the vertex $x_i$ such that there are isometries
  \begin{align*} \alpha_i & : L' \to L_{i-1}, ~~{\rm with}~~ \alpha_i(x_i) = x_{i-1}, \end{align*}
and 
   \begin{align*} \beta_i & : E' \to E_{i-1}, ~~{\rm with}~~ \beta_i(x_i) = x_{i-1}. \end{align*}
Set $L_\delta'$ and $E_\delta'$ to denote the subgraphs of $\cF_2$ containing the vertices $x_ib^{\delta 2^{i-1}}$ and $x_i a^{\delta 2^{i-1}}$ respectively, such that there are isometries
   \begin{align*} \alpha_i^\delta & : L_\delta' \to L', ~~{\rm with}~~ \alpha_i^\delta(x_ib^{\delta 2^{i-1}}) = x_{i}, ~ \delta \in \{-1,1\} ,\end{align*}
and 
   \begin{align*} \beta_i^\delta & : E_\delta' \to E', ~~{\rm with}~~ \beta_i^\delta(x_ia^{\delta 2^{i-1}}) = x_{i}, ~ \delta \in \{-1,1\}. \end{align*}
Let $E_i = E' \cup E_1' \cup E_{-1}' \cup L_1' \cup L_{-1}'$, $L_i = E' \cup L_1' \cup L_{-1}'$, and 
  $$P(\alpha) = \bigcup_{i \in \mN} E_i.$$
We always have $x_0 \in E_i$. If the sequence $\alpha \in \cQ$ is eventually periodic with period $2$, then for large $i$ the complement of $E_i$ in $P(\alpha)$ has $2$ unbounded connected components, and $P(\alpha)$ has $2$ ends. There are also two sequences which are periodic with period $1$, that is, $\alpha = (b)$ and $\alpha = (b^{-1})$. In this case the complement of $E_i$ has one unbounded connected component, and $P(\alpha)$ has $1$ end.

We show that $P(\cQ)$ is a section of $\overline{\cR(K)} \setminus \cR(K)$. Let $V_N = D_X(K,2^{N-1}) \cap X$, that is, 
   $$ (T,e) \in V_N ~ \textrm{if and only if there is an isomorphism} ~ D_K(e,2^{N-1}) \to D_T(e,2^{N-1}).$$
Then $V_1$ hits every orbit in $\overline{\cR(K)}$. Consider $V_2 \subset V_1$. We show that if $(T,e) \in V_2  \setminus \cR(K)$, there exists $\alpha = (\alpha_0,\ldots,\alpha_n(aa^{-1}))$ such that $P(\alpha) \in \cR(T)$.

Indeed, if $(T,e) \in V_2 \setminus \cR(K)$ there is necessarily an isometry $D_{K}(a^{m},2) \to D_{T}(e,2)$, $m \in \mZ$, and $e$ is a point of attachment of a vertical decoration $F$. Since $(T,e) \notin \cR(K)$, by Lemma \ref{lemma:infinitedecorations} $F$ is of finite length $\ell= 2^i-1$, $i \in \mN$. Let $x_0$ be the vertex in $V(F)$ such that $d_T(e,x_0) = \ell$, and either $x_0 = t(w_b)$ or $x_0 = s(w_b)$, where $w_b$ is an edge adjacent to $x_0$. For definitiveness assume that $t(w_b) = x_0$. Then for $0 < k \leq i$ set  
 $$x_k = x_{k-1} \cdot b^{-2^{k-1}}, ~\alpha_{k-1} = b^{-1}.$$
Then $x_i = e$. For $k \geq i$, implement the following inductive procedure, which uses the fact that $(T,e) \in \overline{\cR(K)}$, and so every pattern in $T$ must be replicated in $K$. The boundary of the ball $B_{T}(x_{k},2^{k-1})$ has two adjacent $a$-edges $v_1$ and $v_2$, such that $t(v_1),s(v_2) \subset \partial B_T(x_k,2^{k-1})$ and $s(v_1),t(v_2) \notin B_T(x_k,2^{k-1})$. Let $F_1$ and $F_2$ be vertical decorations attached to $p_1 = s(v_1)$ and $p_2 = t(v_2)$ respectively. Then one of them must have length $2^{k}-1$, and another one must have length $2^{k+1}-1$. If $p_1$ is the vertex with decorations of length $2^{k} - 1$, set $x_{k+1} = p_1$ and $\alpha_{k} = a^{-1}$. Otherwise set $x_{k+1} = p_2$ and $\alpha_{k} = a$. In this way we obtain $\alpha \in \cQ$ such that $P(Q) = (T,e) \cdot b^{2^i-1}$.

Now suppose $(T,e) \in V_1$ and $\cR(T) \cap V_2 = \emptyset$. Let $\{g_s\} \in V(K)$ be a sequence of vertices such that for every $r>0$ there is $s_r>0$ such that for every $s>s_r$ there is an isometry 
   $$\alpha_s^r:D_{K}(g_s,2^r-1) \to D_T(e,2^r-1).$$ 
By choosing a bigger $s_r$, if necessary, we can assume that $D_{K}(g_s,2^r-1) \cap A_0 = \emptyset$, where $A_0$ is the axis of $a$-edges, which means that $D_K(g_s,2^r-1)$ is contained in a vertical decoration $F_s$ of $K$. Then there are two situations: either there exists $r_1>0$ and $s_1 >0$ such that $\partial D_K(g_{s_1},2^{r_1} - 1)$ consists of a single point, and then the same is true for $\partial D_K(g_s,2^r-1)$ for all $s > s_1$ and $r>r_1$, or such an $r_1$ does not exist. In the first case the same condition is true for every $r>r_1$ and $s>s_1$, and it follows that either $(T,e) =P((b^{-1})) \cdot b^{-m}$ or to $P((b)) \cdot b^{m}$ for some $m \in \mZ$. In the second case $(T,e) = P((bb^{-1}))$.
\endproof

Proposition \ref{prop-equicts} completes the proof of Theorem \ref{thm-examples1recurrent}. Recall \cite{Hurder2010} that a pseudogroup $\fG$ of local homeomorphisms of $X$ is \emph{equicontinuous} if for every $\epsilon>0$ there exists $\delta>0$ such that for every $\gamma \in \fG$ and every $x,y \in {\rm dom}(\gamma)$ with $d_X(x,y) < \delta$ we have $d_X(\gamma(x),\gamma(y)) < \epsilon$.

\begin{prop}\label{prop-equicts}
The restriction $\fG|_{V_2}$ is equicontinuous. Every leaf in a graph matchbox manifold $\cM_K$ is at level at most $2$. A leaf $L_T$ is at level $0$ or $1$ if and only if $\cR(T) \cap V_2 = \emptyset$, and there is only a finite number of such leaves.
\end{prop}

\proof Let $L_T$  be a leaf such that $\cR(T) \cap V_2 = \emptyset$. Then by the proof of Lemma \ref{leavesfinitelevel} either $\cR(T) \owns P((bb^{-1}))$, or $\cR(T) \owns P((b))$, or $\cR(T) \owns P((b^{-1}))$. In the first case $L_T$ is a compact leaf and so is at level $0$. In the second and the third case an argument similar to those in Section \ref{sec-totprop} shows that $L_T$ is a one-ended leaf with a single compact leaf in its closure. Such a leaf is at level $1$.

We now show that if $L_T$ is such that $(T,e) \in V_2$, then its orbit under the restricted pseudogroup $\cG_2 = \cG|_{V_2}$ is dense in $V_2$. If that is true, then $L_T$ is dense in $\cM_{K}$ and so is at level $2$. 

To prove this statement we show that $\fG_2$ is equicontinuous. Then by standard arguments, making use of the fact that for an equicontinuous pseudogroup one has uniform control over the domains of maps (see, for example, \cite{ClarkHurder2010b,ALC2009}) the $\fG_2$-orbits of points in $V_2$ are dense in $V_2$.

First we notice that $V_2$ is invariant under the action of the sub-pseudogroup $\cA_4 = \langle \gamma_{a^4} \rangle$, and
   \begin{align}\label{eq-v2st}V_2 \subset {\rm dom}(\gamma_g), ~ \gamma_g \in \cA_4.\end{align} 
Indeed, let $g = a^{4n}$ and $(T,e) \in V_2$. Then $(T,e) \in \cR(P(\alpha_0 \ldots \alpha_n (a a^{-1})))$ and there exists an isomorphism
  $$\alpha: D_{K}(e,2) \to D_T(e,2).$$ 
It follows that $e$ lies on the bi-infinite line of $a$-edges and the action of $g$ on $(T,e)$ is defined. Since $(T,e) \in \overline{\cR(K)}$ there exists $(K,g_n)$ such that there is an isomorphism
  $$\alpha:D_T(e,4(n+1)) \to D_K(g_n,4(n+1)) \in V_2,$$ 
and which implies that $g_n = a^{4m}$. It follows that $(T,e) \cdot g \in V_2$. Next, we claim that $\cG_2 = \cA_4|_{V_2}$. Let $(T,e) \in V_2$ and $g \in F_2$ such that $(T,e) \cdot g \in V_2$. Since $(T,e) \in \overline{\cR(K)}$, then there exists $a^{4 m}$, $m \in \mZ$, such that there is an isomorphism $D_{K}(a^{4 m},\ell(g) + 3) \to D_{T}(e,\ell(g) + 3)$, and so $(K, e) \cdot a^{ 4 m}g$ is defined. Moreover, necessarily 
 \begin{align*} (K, e) \cdot a^{ 4 m}g & = (K, e) \cdot a^{ 4 n}.\end{align*}
Since $K$ is simply connected and translations by $b$ take $(K,e) \cdot a^{ 4 m}$ out of $V_2$, then $g = a^{4(n-m)}$. 

We now show that $\cG_2$ is equicontinuous. Let $(T,e) \in V_2$, and let 
  \begin{equation}2^{i-1}-1 \leq N < 2^i-1, ~ R_N = 2*2^i-1.\end{equation}
Suppose $d_X(T,T') < e^{-R_N}$. We claim that for any $\gamma \in \cG_2$ we have $d_X(\gamma(T),\gamma(T')) < e^{-N}$. Indeed, to see that notice that for $0 \leq 4k < 2^i$ we have an isomorphism 
  $$D_{T \cdot a^{4k}}(e,2^i-1) \to D_T(a^{4k},2^i-1) \subset D_T(e,R_N),$$ 
and therefore there is an isomorphism 
  $$D_{T \cdot a^{4k}}(e,2^i-1) \to D_{T' \cdot a^{4k}}(e,2^i-1).$$ 
Next notice that by construction for every $T \in V_2$ a ball $D_T(e,2^i-1)$ is invariant under the action of $a^{2^i}$. Therefore, for every $k' = 2^i+4k$, $0 \leq 4k < 2^i$ we have an isomorphism 
  $$D_{T \cdot a^{4k}}(e,2^i-1) \to D_{T \cdot a^{k'}}(e,2^i-1),$$ 
and similarly for $T'$. The statement follows. 
\endproof

\section{Some further properties of graph matchbox manifolds}\label{further}

\subsection{Growth of leaves at finite levels}\label{sec-growth}
 
In this section we prove Theorem \ref{thm-examples1growth}, which states that $\fM_n$ contains a leaf of a linear growth, and a leaf of exponential growth. It is enough to construct such examples in $\fM_2$. Let $G_0 = \{a,b\}$. 

We first recall some definitions \cite{CandelConlon2000}. If $(T',e) \in \cR(T) $ define the `plaque distance' by 
  \begin{align*} PD(T',T) & = \min \{ \ell_w(g) ~|~ \gamma_g(T,e) = (T',e) \}, \end{align*}
where $\ell_w(g)$ is the length of $g$ in the word metric on $F_2$, and otherwise $PD(T,T') = \infty$. Then $PD(T,T') = d_{\G_T}(T,T')$. For $k \geq 0$ let
  \begin{align*} \G_k(\cR(T)) & = \{ (T',e) \in \cR(T)  ~|~ PD(T',T) \leq k  \} ,\end{align*}
and define the growth function of $\cR(T)$ by
  \begin{align*} H_{\cR(T)} & : \mZ^+ \to \mR^+ : k \mapsto {\rm card}(\G_k(\cR(T))). \end{align*}
Recall from \cite[Proposition 12.2.35]{CandelConlon2000} that the growth function $H_{\cR(T)}$ can be considered to be the growth function of the leaf $L_T$.

\proof \emph{(of Proposition \ref{thm-examples1growth})}. Let $F_1$ be a subgraph of $\cF_2$ with the set of vertices   
$$V(F_1) = \{ a^n ~|~ n \in \mN \cup\{0\}~\}.$$ 
Then $\cM_{F_1}$ is at level $1$, and $L_{F_1}$ has linear growth, that is,
   \begin{align*} H_{\cR(F_1)}(k) & = k+1. \end{align*}
Let $F_2$ be a subgraph of $\cF_2$ with the set of vertices 
  $$V(F_2) =  \{e\} \cup\{ ag \in F_2 ~| ~ ag ~\textrm{ is a reduced word} ~\}.$$ 
By an argument similar to that in Lemma \ref{leaflevel1} the matchbox manifold $\cM_{F_2}$ consists of two leaves, namely $L_{F_2}$ and a genus two surface $L_{\cF_2}$. Thus $\cM_{F_2}$ is at level $1$. We also have
  \begin{align*} H_{\cR(F)}(k) & = \frac{1}{2} (1 + 3^{k}), \end{align*}
so $L_{F_2}$ has exponential growth. An example of a leaf with polynomial growth which is not totally proper is given in Theorem \ref{thm-examples1recurrent}.
\endproof

\subsection{Pseudogroup dynamics of matchbox manifolds}\label{sec-pseudgroupdynamics}

In this section we study pseudogroup dynamics of graph matchbox manifolds in $\fM_n$ and prove Theorem \ref{thm-dynamicsvsrecurrency}. We first recall a definition of a foliation with expansive dynamics.

\begin{defn}\cite[Section 3]{Hurder2010}
Let $(M,\cF)$ be a foliated space with a foliated atlas 
  $$\cU = \{\varphi_i: U_i \to [-1,1]^n \times \fX_i\},$$ 
and let $\fG$ be the holonomy pseudogroup associated to $\cU$. The dynamics of $\cF$ is $\epsilon$-expansive, or $\fG$ is $\epsilon$-expansive, if there exists $\epsilon>0$ so that for all $w \ne w' \in \fX_i$ with $d_\fX(w,w') < \epsilon$ there exists a holonomy homeomorphism $h \in \fG$ with $w,w' \in {\rm dom} (h)$ such that $d_{\fX}(h(w),h(w')) \geq \epsilon$.
\end{defn}

Proof of statement $(1)$ of Theorem \ref{thm-dynamicsvsrecurrency} is contained in the following lemma.

\begin{lemma}\label{thm:minimalgroup}
Let $L_T \subset \fM_n$ be a non-compact leaf and $\cM_T=\overline{L_T}$ be a graph matchbox manifold. Then for every $0<\epsilon<e^{-2}$ the foliation of $\cM_T$ has $\epsilon$-expansive dynamics. If $L_T$ is compact, then $\cM_T$ is equicontinuous.
\end{lemma}

\proof Fix $\epsilon < e^{-2}$ and recall that the pseudogroup $\fG$ restricted to $X$ is finitely generated with the set of generators $\fG^0$. Suppose there exists $(T,e) \ne (T',e) \in X$ such that for all $\gamma_g \in \fG$ with $(T,e), (T',e) \in {\rm dom} (\gamma_g)$ and $d_X\left((T,e),(T',e)\right) < \epsilon$ we have $d_X(\gamma_g(T,e),\gamma_g(T',e)) < \epsilon$. Notice that since $\epsilon < e^{-2}$ and $d_X((T,e),(T',e)) < \epsilon$, then for all $\gamma_i \in \fG^0$ $(T,e) \in {\rm dom} (\gamma_i)$ if and only if $(T',e) \in {\rm dom} (\gamma_i)$. We show that in fact this implies that $(T,e) \in {\rm dom} (\gamma_g)$ if and only if $(T',e) \in {\rm dom}(\gamma_g)$ for all $\gamma_g \in \fG$.

We have $\gamma_g = \gamma_{i_1} \circ \cdots \gamma_{i_k}$. The proof is by induction on $\ell$, $1 \leq \ell \leq k$. Suppose $(T,e),(T',e) \in {\rm dom} (\gamma_{i_1} \circ \cdots \circ \gamma_{i_\ell})$. Set 
 $$(T_\ell, e) = \gamma_{i_1} \circ \cdots \circ \gamma_{i_\ell}(T,e), ~ {\rm and} ~ (T_\ell', e) = \gamma_{i_1} \circ \cdots \circ \gamma_{i_\ell}(T',e).$$ 
An assumption that $d_X((T_\ell,e),(T_\ell',e)) < \epsilon$ means that there is an isometry 
   $$\alpha_\ell:D_{T_\ell}(e,2) \to D_{T_\ell'}(e,2),$$  
and $(T_\ell,e) \in {\rm dom} (\gamma_{i_{\ell+1}})$ if and only if $(T_\ell',e) \in {\rm dom}(\gamma_{i_{\ell+1}})$. Therefore, $(T,e)\in {\rm dom} (\gamma_{i_1} \circ \cdots \circ \gamma_{i_{\ell+1}})$ if and only if $(T',e) \in {\rm dom} (\gamma_{i_1} \circ \cdots \circ \gamma_{i_{\ell+1}})$.

We now show that this implies $(T,e) = (T',e)$. Let $R \in \mN$. We want to show that there is an isomorphism $\alpha_R: D_T(e,R) \to D_{T'}(e,R)$. The proof is by induction on $r$, $1 \leq r \leq R$. Suppose there is an isomorphism   $$\alpha_r: D_T(e,r) \to D_{T'}(e,r)$$ 
and let $\gamma_{g_1},\ldots,\gamma_{g_{n_r}}$ be a collection of homeomorphisms such that $\ell(g_k) = r$ and $(T,e),(T',e) \in {\rm dom} (\gamma_{g_k})$. Notice that 
 $$D_T(e,r+1) \subset D_T(e,r) \cup \left( \bigcup_{1 \leq k \leq n_r} D_T(g_{k},2) \right)$$ 
and similarly for $T'$, so if for every $1 \leq k \leq n_r$ we have $D_T(g_k,2) = D_{T'}(g_k,2)$, then we are done. But the latter condition follows from the fact that 
 \begin{align*} d_X \left(\gamma_{g_k}(T,e),\gamma_{g_k}(T',e) \right) < \epsilon < e^{-2}. \end{align*}
Now if $L_T$ is non-compact, then it contains at least one other leaf $L_S$, $(S,e) \in \overline{\cR(T)} \setminus \cR(T)$. Let $(T',e) \in \cR(T)$ such that $d_X((T',e),(S,e)) < \epsilon$. By the argument above $(S,e) = (T',e) \in \cR(T)$, which contradicts the assumption. If $L_T$ is compact, then $\cR(T)$ is a finite set, and so there exists $\epsilon>0$ such that for any $(T',e) \in \cR(T)$ the $\epsilon$-neighborhood of $(T',e)$ in $\overline{\cR(T)} = \cR(T)$ contains exactly one point. Then $\cM_T$ is trivially equicontinuous. 
\endproof

We now restate and prove corollary \ref{no-solenoid}.

\begin{cor}
The foliated space $\fM_n$ does not contain a weak solenoid.
\end{cor}

\proof A weak solenoid must contain a non-compact leaf and have an equicontinuous pseudogroup, which is not possible by Lemma \ref{thm:minimalgroup}. \endproof

The proof of statement $(2)$ of Theorem \ref{thm-dynamicsvsrecurrency} is given by the following lemma. 

We say that the pseudogroup $\fG$ restricted to a transversal $X$ has a non-trivial equicontinuous subsystem if and only if there exists a clopen subset $C \subset X$ such that $\fG_C = \fG|_C$ is equicontinuous.

\begin{lemma}
Let $\fG_T$ be the holonomy pseudogroup associated to the transversal $\overline{\cR(T)}$ of a graph matchbox manifold $\cM_T$. If $\fG_T$ has a non-trivial equicontinuous subsystem, then $L_T$ is a recurrent leaf. The converse is false, i.e. there exists a matchbox manifold $\cM_{T'}$ where $L_{T'}$ is recurrent and $\fG_{T'}$ has no nontrivial equicontinuous subsystem.
\end{lemma}

\proof Let $C \subset \overline{\cR(T)}$ be a clopen neighborhood with equicontinuous restricted pseudogroup $\fG_C$. By standard methods one can show that equicontinuity of $\fG_C$ implies that there is a uniform lower bound on the domains of holonomy homeomorphisms, that is, there exists $\eta>0$ such that if $x \in D(h)$, $h \in \fG_C$, then $D_X(T,\eta) \subset D(h)$ (see, for instance, \cite{ClarkHurder2010b,ALC2009}). It follows that orbits of points in $C$ under $\fG_C$ are dense in $C$, and so $\cR(T)$ accumulates on itself.

Conversely, let $\cM_{T'}$ be minimal, and suppose there is a clopen neighborhood $C \subset \overline{\cR(T')}$ with equicontinuous restricted pseudogroup $\fG_C$. Then $\fG$ must be equicontinuous which contradicts Lemma \ref{thm:minimalgroup}.
\endproof



\begin{thebibliography}{10}

\bibitem{AO1995}
{J.~Aarts and L.~Oversteegen},
\newblock {\it Matchbox manifolds},
\newblock In {\bf Continua ({C}incinnati, {OH}, 1994)},
\newblock {Lecture Notes in Pure and Appl. Math., Vol. 170},
\newblock {Dekker, New York}, 1995, pages 3--14..

\bibitem{AM1988}
{J.M.~Aarts and M.~Martens},
\newblock {\it Flows on one-dimensional spaces},
\newblock {\bf Fund. Math.}, 131:39--58, 1988.

\bibitem{ALM2008}
{F.~Alcalde~Cuesta, A.~Lozano~Rojo and M.~Macho~Stadler},
\newblock {\it Dynamique transverse de la lamination de {G}hys-{K}enyon},
\newblock {\bf Ast\'erisque}, 323:1--16, 2009.

\bibitem{ALC2009}
\newblock {J.~A.~{\'A}lvarez L{\'o}pez and A.~Candel},
 \newblock {\it Equicontinuous foliated spaces},
\newblock {\bf Math. Z.}, {263(4)}, 725--774, 2009.

\bibitem{BG2003}
{R.~Benedetti and J.-M.~Gambaudo},
\newblock {\it On the dynamics of $\mG$-solenoids. {A}pplications to {D}elone sets},
\newblock {\bf Ergodic Theory Dyn. Syst.}, 23:673--691, 2003.

\bibitem{Blanc2001}
{ E.~Blanc},
\newblock {\bf Propri\'et\'es g\'en\'eriques des laminations},
\newblock  {PhD Thesis}, {Universit\'e de Claude Bernard-Lyon 1},  {Lyon}, 2001.

\bibitem{Burago}
{D.~Burago, Yu. ~Burago and S.Ivanov},
\newblock {\bf A Course in Metric Geometry},
\newblock GSM 33, American Mathematical Society, Providence, Rhode Island 2001.

\bibitem{CandelConlon2000}
{A.~Candel and L.~Conlon},
\newblock {\bf Foliations I},
\newblock Amer. Math. Soc., Providence, RI, 2000.

\bibitem{CandelConlon2003}
{A.~Candel and L.~Conlon},
\newblock {\bf Foliations II},
\newblock Amer. Math. Soc., Providence, RI, 2003.

\bibitem{CC1}
  {J. ~Cantwell and L. ~Conlon},
     \newblock{\it Poincar\'e-{B}endixson theory for leaves of codimension one},
 \newblock {\bf Trans. Amer. Math. Soc.}, 265(1):181--209, 1981.

\bibitem{CC3}
  {J. ~Cantwell and L. ~Conlon},
     \newblock{\it Growth of leaves},
 \newblock {\bf Comment. Math. Helv.}, 53:93--111, 1978.

  
\bibitem{ClarkHurder2010a}
{A.~Clark and S.~Hurder},
\newblock {\it Embedding solenoids in foliations},
\newblock {\bf Topology Appl.}, 158:1249--1270, 2011.

\bibitem{ClarkHurder2010b}
{A.~Clark and S.~Hurder},
\newblock {\it Homogeneous matchbox manifolds},
\newblock {to appear in {\bf Transactions AMS}}, version 2011.

\bibitem{EMT1977}
{D.B.A.}~Epstein, {K.C.}~Millet, and {D.}~Tischler,
\newblock {\it Leaves without holonomy},
\newblock {\bf Jour. London Math. Soc.}, 16:548--552, 1977.

\bibitem{FO2002}
{R.~Fokkink and L.~Oversteegen},
\newblock {\it Homogeneous weak solenoids},
\newblock {\bf Trans. Amer. Math. Soc.}, 354:3743--3755, 2002.

\bibitem{Ghys1999}
{\'{E}~Ghys},
\newblock {\it Laminations par surfaces de Riemann},
\newblock In {\bf Dynamique et G\'{e}om\'{e}trie Complexes},
\newblock {Panoramas \& Synth\`{e}ses}, 8:49--95, 1999.

\bibitem{Hector}
{G. ~Hector},
\newblock {\it Architecture des feuilletages de class $C^2$},
\newblock {Asterisque, 107-108, Soci\'ete Math\'ematique de France}, 1983: 243-258.

\bibitem{Hector2}
{G. ~Hector},
\newblock {\it Feuilletages en cylindres}, 
\newblock {\bf Lecture Notes in Math.}, 597, Springer-Verlag, 1977, pp. 252-270.

\bibitem{Hector3}
{G. ~Hector},
\newblock {\it Quelques examples des feuilletages - Esp\'eces rares}, 
\newblock {\bf Ann. Inst. Fourier}, 26(1), 1976, pp. 239-264.

\bibitem{Hurder2010}
{S.~Hurder},
\newblock {\it Lectures on foliation dynamics: Barcelona 2010}, preprint 2011.


\bibitem{Kur}
{K.~Kuratowski},
\newblock {\bf Topology I},
\newblock {Academic Press Inc.}, 1966, 547pp.


\bibitem{LozanoSp}
{A.~ Lozano Rojo},
\newblock {\it Foliated spaces defined by graphs},
\newblock Rev. Semin. Iberoam. Mat.. 3(4):21--38, 2007.


\bibitem{Lozano2009}
{A.~ Lozano Rojo},
\newblock {\it An example of a non-uniquely ergodic lamination},
\newblock preprint 2009.

\bibitem{Salhi2003}
{H.~Marzougui and E.~Salhi},
\newblock {\it Structure of foliations of codimension greater than one},
\newblock {\bf Comment. Math. Helv.}, 78(4):722--730, 2003.

\bibitem{McCord1965}
{C.~McCord},
\newblock {\it Inverse limit sequences with covering maps},
\newblock {\bf Trans. Amer. Math. Soc.}, 114:197--209, 1965.

\bibitem{MS2006}
{C.C.~Moore and C.~Schochet},
\newblock {\bf Global analysis on foliated spaces},
\newblock {With appendices by S. Hurder, Moore, Schochet and Robert J.
              Zimmer}
\newblock Math. Sci. Res. Inst. Publ. vol. 9,
\newblock Springer-Verlag, New York, 1988, vi+337 pp.


\bibitem{Nish1977}
{T.~Nishimori},
\newblock{\it Behaviour of leaves of codimension-one foliations},
\newblock T\^ohoku Math. J., 29(2):255--273, 2977.

\bibitem{Nish1979} 
{T.~Nishimori},
\newblock{\it Ends of leaves of codimension-one foliations},
\newblock{T\^ohoku Math. J. (2), 1979(1):1--22, 1979.}

\bibitem{FrankSadun2009}
{N. ~ Priebe Frank and L.~Sadun}
\newblock {\it Topology of some tiling spaces without finite local
              complexity},
\newblock {\bf Discrete Contin. Dyn. Syst.}, 23(3):847--865, 2009.

\bibitem{Salhi1982}
{E.~Salhi},
\newblock {\it Sur les ensembles minimaux locaux},
\newblock {\bf C. R. Acad. Sci. Paris S\'er. I Math.}, 295(12):691--694, 1982.

\bibitem{Salhi1985-1}
{E.~Salhi},
\newblock {\it Niveau des feuilles},
\newblock {\bf C. R. Acad. Sci. Paris S\'er. I Math.}, 301(5):219--222, 1985.

\bibitem{Salhi1985-2}
{E.~Salhi},
\newblock {\it Sur un th\'eor\`eme de structure des feuilletages de
              codimension {$1$}},
\newblock {\bf C. R. Acad. Sci. Paris S\'er. I Math.}, 300(18):635--638, 1985.

\bibitem{Sullivan1988}
{D.~Sullivan},
\newblock{Bounds, quadratic differentials, and renormalization conjectures},
\newblock{American {M}athematical {S}ociety centennial publications, {V}ol. {II} ({P}rovidence, {RI}, 1988)},
    {417--466},{1992}

\bibitem{Ts1980}
{N.~Tsuchiya}
\newblock {\it Leaves of finite depth},
\newblock {\bf Japan. J. Math.}, 6(2):343--364, 1980.

\bibitem{Williams1974}
{R.F.~Williams},
\newblock {\it Expanding attractors},
\newblock {\bf Inst. Hautes \'Etudes Sci. Publ. Math.}, 43:169--203, 1974.

\end{thebibliography}
\end{document}